\documentclass[11pt]{article}
\usepackage{amscd}
\usepackage{amsmath}
\usepackage{latexsym}
\usepackage{amsfonts}
\usepackage{amssymb}
\usepackage{amsthm}
\usepackage{graphicx}

\newtheorem{theorem}{Theorem}[section]
\newtheorem{definition}{Definition}[section]
\newtheorem{proposition}{Proposition}[section]
\newtheorem{lemma}{Lemma}[section]
\newtheorem{remark}{Remark}[section]

\makeatletter
\newcommand{\Extend}[5]{\ext@arrow0099{\arrowfill@#1#2#3}{#4}{#5}}
\makeatother

\begin{document}
\author{}
\title{Scattering for the focusing ${\dot H}^{1/2}$-critical
 Hartree equation with  radial data}
\author{{Yanfang Gao$^{1}$, \;  Changxing Miao$^2$\; and \; Guixiang Xu$^2$ }\\
    {\small $^1$ Institute of  Mathematics, Jilin University, Changchun,
    China, 130012}\\
         {\small $^2$ Institute of Applied Physics and Computational Mathematics}\\
         {\small P. O. Box 8009,\ Beijing,\ China,\ 100088}\\
         {\small ( gaoyanfang236@yahoo.com.cn,\; miao\_changxing@iapcm.ac.cn, \; xu\_guixiang@iapcm.ac.cn) }\\
         \date{}
        }
\maketitle

\maketitle

\begin{abstract}
We investigate the focusing  $\dot H^{1/2}$-critical  nonlinear
Schr\"{o}dinger equation (NLS) of Hartree type $i\partial_t u +
\Delta u = -(|\cdot|^{-3} \ast |u|^2)u$ with $\dot H^{1/2}$ radial
data in dimension $d = 5$. It is proved that if the maximal
life-span solution obeys
$\sup_{t}\big\||\nabla|^{\frac{1}{2}}u\big\|_2 <
\frac{\sqrt{6}}{3} \big\||\nabla|^{\frac{1}{2}}Q\big\|_2$, where
$Q$ is the positive radial solution to the elliptic equation with
nonlocal operator (\ref{e14}) which  corresponds to a new
variational structure. Then the solution is global and scatters.
\end{abstract}

 \begin{center}
 \begin{minipage}{120mm}
 \end{minipage}
 \end{center}

\section{Introduction}
\setcounter{section}{1}\setcounter{equation}{0}

Consider the Cauchy problem for the  $\dot H^{1/2}$-critical
Hartree equation
\begin{equation}
i\partial_t u + \Delta u = F(u)
\end{equation}
in $\mathbb R^5$, where $F(u) = -(|\cdot|^{-3} \ast |u|^2)u$, $u$
is a complex-valued function defined on some spacetime slab $I
\times \mathbb R^5$. The Hartree equation arises in the study of
boson stars and other physical phenomena, see, for instance,
\cite{c00}.

The term $\dot H^{1/2}$-critical means that the scaling
\begin{equation}
u_\lambda(t,x) = \lambda^{-2}u(\lambda^{-2}t, \lambda^{-1}x)
\end{equation}
leaves both the equation and the initial data of $\dot H^{1/2}_x$-
norm invariant. By a function $u : I \times \mathbb R^5 \mapsto
\mathbb C$ is a {\it solution} to (1.1), it means that $u \in
C_t^0\dot H^{1/2}_x(K \times \mathbb R^5) \cap L_t^3L_x^{15/4}(K
\times \mathbb R^5)$ for any compact $K \subset I$, and $u$ obeys
the Duhamel formula
\[u(t) = e^{i(t-t_0)\Delta}u(t_0) - i \int_{t_0}^t e^{i(t-t')\Delta} F(u(t'))\,\mathrm{d} t' \]
for all $t,\, t_0 \in I$. We call $I$ the {\it life-span} of $u$.
If $I$ can not be extended strictly larger, we say $I$ is the {\it
maximal life-span} of $u$, and $u$ is a {\it maximal life-span
solution}. If $I = \mathbb R$, then $u$ is global.
\begin{definition}[Blow up]
 Let $u : I \times \mathbb R^d \mapsto \mathbb C$ be a solution to (1.1). Say $u$
blows up forward in time if there exists $t_1 \in I$ such that
\[\|u\|_{L_t^3L_x^{15/4}([t_1, \; \sup I)\times \mathbb R^5)} =\infty \,;\]
and $u$ blows up backward in time if there exists $t_1$ such that
\[ \|u\|_{L_t^3L_x^{15/4}((\inf I, \; t_1]\times \mathbb
R^5)}=\infty.\]
\end{definition}

Throughout the paper, we write
\[\|u\|_{S(I)} := \|u\|_{L_t^3L^{15/4}_x(I \times \mathbb R^5)},
 \quad \|u\|_{X(I)} := \big\||\nabla|^{\frac{1}{2}}u\big\|_{L_t^3L_x^{30/11}(I \times \mathbb R^5)}.\]

The local theory for (1.1) was established by Cazenave and
Weissler \cite{c3}, \cite{c4}. Using a fixed point argument
together with Strichartz's estimates in the framework of Besov
spaces, they constructed local in time solution for arbitrary
initial data. However, due to the critical nature of the equation,
the existence time depends on the profile of the initial data and
not merely on its $\dot H^{1/2}_x$-norm. They also proved the
global existence for small data.
\begin{theorem}[Local theory, \cite{c3}, \cite{c4}]
 Let $u_0 \in \dot H_x^{1/2}(\mathbb R^5), \, t_0 \in
\mathbb R$, there exists a unique maximal life-span solution $u :
I \times \mathbb R^5 \mapsto \mathbb C$ to $(1.1)$ with initial
data $u(t_0) = u_0$. This solution also has the following
properties: \vspace{-8pt}
\begin{itemize}
\item(Local existence) $I$ is an open neighborhood of $t_0$.
\item(Blow up criterion) If $\sup I$ is finite, then $u$ blows up
forward in time; if $\inf I$ is finite, then $u$ blows up backward
in time. \item(Scattering) If $\sup I = +\infty$, and $u$ does not
blow up forward in time, then $u$ scatters forward in time, that
is, there exists a unique $u_+ \in \dot H^{1/2}_x(\mathbb R^5)$
such that \begin{equation} \lim_{t \to +\infty} \|u(t) -
e^{it\Delta}u_+\|_{\dot H^{1/2}_x(\mathbb R^5)} = 0.
\end{equation}
Conversely, given $u_+ \in \dot H^{1/2}_x(\mathbb R^5)$, there
exists a unique solution to $(1.1)$ in a neighborhood of infinity
such that $(1.3)$ holds. \item(Small data scattering) If
$\big\||\nabla|^{\frac{1}{2}}u_0\big\|_2$ is sufficiently small,
then $u$ scatters in both time directions. Indeed,
$\|u\|_{S(\mathbb R)} \lesssim
\big\||\nabla|^{\frac{1}{2}}u_0\big\|_2$. \item(Radial symmetry)
If $u_0$ is radially symmetric, then $u$ remains radially
symmetric for all time.
\end{itemize}
\end{theorem}
From Theorem 1.1,  a solution to (1.1) with small data must be
scattering. However, the result is unknown  for arbitrary data,
even in the defocusing case. In \cite{c11}, Kenig and Merle proved
for the defocusing cubic NLS that the solution is global and
scatters if it remains uniformly bounded in $\dot H^{1/2}_x$ on
its maximal life-span. The assumption that the solution is
uniformly bounded in $\dot H^{1/2}_x$ plays a role of the missing
conservation law. The argument presented there applies to the
corresponding defocusing Hartree equation without difficulty. As
to the focusing case, there has been no result on the line of
scattering, neither NLS nor of  Hartree type. Our primary goal in
this paper is to establish scattering result for the focusing
Hartree equation, and we believe that the argument can be adapted
to the focusing NLS.

For the Cauchy problem of $(1.1)$, there is a stationary solution
 $e^{it}\bar Q$  that is global but blows up both
forward and backward.  Here $\bar Q$ is  the unique positive
radial Schwartz solution to
\[ \Delta \bar Q + (|\cdot|^{-3} \ast |\bar Q|^2)\bar Q = \bar Q.\]

In the focusing energy/mass critical case, the corresponding
stationary solution/ground state  play the role of an obstruction
to the global well-posedness and scattering. Indeed, the global
existence follows  so long as the kinetic energy/mass of the
initial data is strictly less than that of the stationary
solution/ground state. In \cite{c23}, Li-Zhang classify the
minimal blowup solutions of the focusing mass-critical Hartree
equation.  However, wether the solution $u $ to $(1.1)$ on its
maximal life-span with $\|u\|_{L_t^\infty\dot H_x^{1/2}} < \|\bar
Q\|_{\dot H^{1/2}_x}$ implies global existence is still open. In
this paper  we  will introduce a new elliptic equation:
\begin{equation}\label{e14}
\Delta Q + (|\cdot|^{-3} \ast  |Q|^2)Q = (-\Delta)^{1/2} Q,
\end{equation}
which corresponds to a new variational structure, and prove that
if the solution $u$ to $(1.1)$ satisfies $\|u\|_{L_t^\infty\dot
H_x^{1/2}} < \frac{\sqrt{6}}{3}\| Q\|_{\dot H^{1/2}_x}$, then the
solution is global and scatters.

Solutions to critical NLS and of Hartree type have been
intensively studied, especially those of energy critical
equations. Scattering results for the defocusing energy-critical
equations have been completely established. These were
accomplished by Bourgain \cite{c2}, Grillakis \cite{c8}, Tao
\cite{c24}, Colliander-Keel-Staffilani-Takaoka-Tao \cite{c5},
Ryckman-Visan \cite{c21}, and Visan \cite{c27}, Miao-Xu-Zhao
\cite{c19}. As will be discussed later, the focusing
energy-critical NLS theory has also been well established by
Kenig-Merle and Killip-Visan, except for dimensions 3 and 4. For
the focusing Hartree, it was proved by Li-Miao-Zhang \cite{c16},
and Miao-Xu-Zhao\cite{c24}.

Another kind  of critical NLS and of Hartree type which receives
lots of attention is the mass-critical one.  Results in earlier
work which is devoted to global well-posedness were usually
obtained under the assumption of the $H^1_x$ initial data. See,
e.g., \cite{c3}, \cite{c28}. In \cite{c28}, Weinstein first
observed the role of the ground state for the focusing
mass-critical NLS despite finite energy. As far as $L_x^2$ initial
data is concerned, Tao-Visan-Zhang \cite{c25} proved the
scattering results for the defocusing case for large spherically
symmetric data in dimensions three and higher. More recent and
nice work on scattering results for $L_x^2$ data were done by
Killip-Tao-Visan
 \cite{c13}, Killip-Visan-Zhang \cite{c15}, and Miao-Xu-Zhao \cite{c20} with spherical
symmetry assumption.

The recent progress in studying those equations is due to a new
and highly efficient approach based on a concentration compactness
idea to provide a linear profile decomposition. This approach
arises from investigating the defect of compactness for the
Strichatz estimates. Based on a refined Sobolev inequality,
Kerrani \cite{c12} obtained a linear profile decomposition for
solutions of free NLS with $H^1_x$ data.  It was Kenig and Merle
who first introduced Kerrani's linear profile decomposition to
obtain scattering results. They treated the focusing
energy-critical NLS in dimensions 3, 4, 5 in \cite{c10}. Using the
same decomposition, Killip and Visan \cite{c14} dealt with the
focusing energy-critical NLS in dimensions five and higher without
radial assumption. Using the decomposition of \cite{c18},
Tao-Visan-Zhang \cite{c26} made a reduction for failure of
scattering. And by combining the reduction with an in/out
decomposition technique, \cite{c13}, \cite{c15} settled the
scattering problem for the mass-critical NLS  with spherically
symmetric data.

A linear profile decomposition for general $\dot H^s$ data was
proved by Shao \cite{c22}. Unlike Kerrani's approach which is
based on a refined Sobolev inequality, Shao took advantage of the
existing $L_x^2$ linear profile decomposition and the Galilean
transform, and managed to eliminate the frequency parameter from
the decomposition. In this paper, we will use Shao's linear
profile decomposition, and our main result is:
\begin{theorem}
Let $u_0 \in \dot H^{1/2}_x(\mathbb R^5)$, radially symmetric, $t_0
\in \mathbb R$, $I$ is a time interval containing $t_0$. Let $u : I
\times \mathbb R^5 \mapsto \mathbb C$ be a maximal life-span
solution to $(1.1)$. Assume $\sup_{t \in I}
\big\||\nabla|^{\frac{1}{2}}u(t)\big\|_2  < \frac{\sqrt{6}}{3}
\big\||\nabla|^{\frac{1}{2}}Q\big\|_2$. Then $u$ is global and
scatters with
\[ \|u\|_{L_t^3 L_x^{15/4}(\mathbb R \times \mathbb R^5)}^3 =
\int_{\mathbb R} \left(\int_{\mathbb R^5} |u(t,x)|^{15/4}
\,\mathrm{d} x\right)^{4/5} \,\mathrm{d} t < \infty.\]
\end{theorem}

\begin{remark}  It is  an interesting problem
to describe the correspondence between $Q$ and $\bar Q$, and thus
leading to some investigation with the gap.  It is also an
interesting problem that wether the solution blows up so long as
$\sup\limits_{t \in I} \big\||\nabla|^{\frac{1}{2}}u(t)\big\|_2 \geq
\frac{\sqrt{6}}{3} \big\||\nabla|^{\frac{1}{2}}Q\big\|_2$.
\end{remark}

The concentration compactness argument reduces matters to the
study of almost periodic solutions modulo symmetries.

\begin{definition}[Almost periodic modulo scaling]
 Let $u$ be a solution to $(1.1)$ with maximal life-span $I$. Say
$u$ is almost periodic modulo scaling if there exist functions $N
: I \mapsto \mathbb R^+$, $C : \mathbb R^+ \mapsto \mathbb R^+$
such that for all $\eta > 0$, $t \in I$
\[\int_{|x| \geq C(\eta)/{N(t)}} \big||\nabla|^{\frac{1}{2}}u(t,x)\big|^2 \,\mathrm{d} x \leq \eta\]
and
\[\int_{|\xi| \geq C(\eta)N(t)} |\xi||\hat u(t,\xi) |^2 \,\mathrm{d} \xi \leq \eta.\]
We refer to $N(t)$ as the frequency scale function for the
solution, and $C$ the compactness modulus function.
\end{definition}

\begin{remark}
By the Arzela-Ascoli theorem, a family of functions is precompact
in $\dot H^{1/2}_x(\mathbb R^5)$ if and only if it is norm-bounded
and there exists a compactness modulus function $C$ so that
\[\int_{|x| \geq C(\eta)}\big||\nabla|^{\frac{1}{2}}f(x)\big|^2 \,\mathrm{d}x + \int_{|\xi| \geq C(\eta)}|\xi||\hat
f(\xi)|^2 \,\mathrm{d} \xi \leq \eta\] for all functions in the
family and all $\eta >0$. Thus, $u$ is almost periodic modulo
scaling if and only if there exists a compact subset $K$ of $\dot
H^{1/2}_x(\mathbb R^5)$ such that
\[\big\{\, u(t): t \in I\, \big\} \subseteq \big\{\,\lambda^{-2}f(\lambda^{-1}x): \lambda \in (0, +\infty), f \in K \, \big\}.\]
By Sobolev's embedding theorem, any solution $u : I \times \mathbb
R^5 \mapsto \mathbb C$ to $(1.1)$ that is almost periodic modulo
scaling also satisfies
\begin{equation}\label{e200}
\int_{|x| \geq C(\eta)/{N(t)}}|u(t,x)|^{\frac{5}{2}} \,\mathrm{d}
x \leq \eta\end{equation} for all $t \in I$ and all $\eta > 0$.

By the compactness modulo scaling, there also exists a function
$c: \mathbb R^+ \mapsto \mathbb R^+$ such that
\begin{equation}\label{e15}
\int_{|x| \leq c(\eta)/{N(t)}}
\big||\nabla|^{\frac{1}{2}}u(t,x)\big|^2 \,\mathrm{d} x +
\int_{|\xi|\leq c(\eta)N(t)}|\xi||\hat u(t,\xi)|^2 \,\mathrm{d}
\xi \leq \eta\end{equation} for all $t \in I$ and all $\eta > 0$.
\end{remark}
We now present the process of reduction. If Theorem 1.2 failed,
then there must be an almost periodic solution. More precisely, we
have:
\begin{theorem}
Suppose Theorem $1.2$ failed for radially symmetric data. Then
there exists a maximal life-span solution $u : I \times \mathbb
R^5 \mapsto \mathbb C$ to $(1.1)$ with $\sup_t
\big\||\nabla|^{\frac{1}{2}}u\big\|_2 < \frac{\sqrt{6}}{3}
\big\||\nabla|^{\frac{1}{2}}Q\big\|_2$. $u$ is almost periodic
modulo scaling, blows up both forward and backward. Moreover, the
frequency scale function $N(t)$ and the maximal life-span $I$
match one of the following scenarios :

{\rm I}. (Finite-time blowup) Either $|\inf I| < \infty$ or $\sup
I < \infty$.

{\rm II}. (Low-to-high cascade) $I = \mathbb R$,
\[\inf N(t) \geq 1 \quad \textrm{for all} \,\,  t \in \mathbb R,\quad \textrm{and} \quad \limsup_{t \to +\infty} N(t) = +\infty.\]

{\rm III}. (Soliton-like solution) $I =\mathbb R$,  $N(t) \equiv
1$ for all $t \in \mathbb R$.
\end{theorem}

The delicate relationship between the frequency scale function and
the maximal life-span for almost periodic solution was first
discovered by Killip, Tao, and Visan in \cite{c13} for
mass-critical NLS. The argument was adapted to the energy-critical
case in \cite{c14}. This latter argument is directly applicable to
the setting of this paper.

To prove Theorem 1.2, it suffices to preclude the three scenarios
in Theorem 1.3. We adapt ideas  in \cite{c13}, \cite{c14}.
However, when precluding the finite-time blowup, Plancherel's
theorem and Hardy's inequality are not enough to obtain a decay
for the localized mass, especially for large scales, as we are
working in the fractional Sobolev space. To surmount this, we take
advantage of the intrinsic description of fractional derivatives,
estimate the integral formula in cases according to the spatial
scales. Some negative regularity is needed for disproving the rest
two scenarios, and our discussions are somewhat involved due to
the nonlocal nonlinearity and low regularity. We shall make full
use of the frequency localization. For instance, in the proof of
Lemma 6.1, we should firstly use Bernstein's inequality to obtain
a positive gain in estimating the high frequency components and
the medium frequency components, such that the Gronwall's
inequality is applicable. What we would also like to emphasize in
particular is that as the $\dot H^{1/2}$-critical equation enjoys
no conservation law, beside proving the negative regularity, we
have to gain additional regularity of at least 1 order
differentiability, which means that the soliton-like solution has
conserved energy;  and thus allows us to apply virial-type
argument to disprove it. We also obtain the local spacetime bounds
in terms of the frequency scale function for all $\dot
H^{1/2}$-admissible pairs and of those $L^2$-admissible pairs
$(q,\,r)$ with $q \geq 3$, $r \leq 30/11$.

The following lemma plays an important role in proving the
negative and additional regularity. See \cite{c26} for a proof.
\begin{lemma}
Let $u$ be an almost periodic solution to $(1.1)$ on its maximal
life-span $I$. Then, for all $t \in I$
\begin{eqnarray}\label{e16}
u(t) & =  & \lim\limits_{T \nearrow  \sup I} i \int_t^T
e^{i(t-t')\Delta}F(u(t'))\,\mathrm{d} t' \nonumber\\ & = & -\lim\limits_{T
\searrow \inf I} \int_T^t e^{i(t-t')\Delta}F(u(t'))\,\mathrm{d} t'
\end{eqnarray}
as weak limits in $\dot H^{1/2}_x$.
\end{lemma}

The rest of paper is organized as follows. In Section 2, we list
out some notations and known results that we use repeatedly in the
paper. In Section 3, the sharp constant for a
Hardy-Littlewood-Sobolev type inequality is obtained, and a
sufficient condition for global existence of $(1.1)$ with finite
energy initial data is given. In Section 4, we first prove a
Palais-Smale condition modulo scaling, and then Theorem 1.3. In
Section 5, we preclude the finite-time blowup scenario. In Section
6, we prove the negative regularity for global case. In Section 7,
we disprove the low-to-high cascade. In Section 8, we prove an
additional regularity for the soliton-like solution. In Section 9,
we preclude the soliton-like solution. In Section 10, we prove
Proposition 1.1.

\section{Preliminaries}
\setcounter{equation}{0}

\subsection{Notations}

For any spacetime slab $I \times \mathbb R^5$, we use
$L_t^qL_x^r(I \times \mathbb R^d)$ to denote the Banach space with
norm
\[\|u\|_{L_t^qL_x^r} := \left(\int_I \left(\int_{\mathbb R^d}|u(t,x)|^r \,\mathrm{d} x\right)^{q/r} \,\mathrm{d} t\right)^{1/q},\]
with the usual modifications when $q $ or $r$ are infinity. When
$q =r$ we abbreviate $L_t^q L_x^r$ as $L_{t,x}^q$.

We use the `Japanese bracket' convention $\langle x \rangle := (1
+ |x|^2)^{1/2}$.

We use $X \lesssim Y$ or $Y \gtrsim X$ whenever $X \leq CY$ for
some constant $C > 0$. If $C$ depends on some parameters, we will
indicate this with subscripts; for example, $X\lesssim_u Y$ denote
 the assertion that $X \leq C_u Y$ for some $C_u$ depending on
 $u$. We denote by $X{\pm}$ any quantity of the form $X \pm
 \varepsilon$ for any $\varepsilon >0$.
we define the Fourier transform on $\mathbb R^d$ by
\[\hat f(\xi) := (2\pi)^{-\frac{d}{2}}\int_{\mathbb R^d}e^{-ix\cdot \xi}f(x) \,\mathrm{d} x .\]
For $s \in \mathbb R$, we define the fractional
differential/integral operators
\[\widehat{|\nabla|^s f}(\xi) := |\xi|^s \hat f(\xi)\]
and the homogeneous Sobolev norm
\[\|f\|_{\dot H^s_x(\mathbb R^d)} := \big\||\nabla|^s f\big\|_{L_x^2(\mathbb R^d)}.\]

The next following lemma is a form of Gronwall's inequality that
we will use to handle some bootstrap argument below.
\begin{lemma}[Gronwall's inequality]
Given $\gamma > 0$, $0 < \eta < \frac{1}{2}(1-2^{-\gamma})$ and
$\{b_k\} \in  l^\infty(\mathbb Z^+) $. Let $\{x_k\} \in
l^\infty(\mathbb Z^+) $ be a non-negative sequence obeying
\[
x_k \leq b_k + \eta\sum_{l =0}^\infty 2^{-\gamma|k - l|}x_l  \quad
\textrm{for all} \,\, k \geq 0. \] Then
\begin{equation}
x_k \lesssim \sum_{l=0}^\infty r^{|k-l|}b_l \quad \textrm{for all}
\,\, k\geq 0
\end{equation}
for some $r=r(\eta) \in (2^{-\gamma}, 1)$. Moreover, $r \downarrow
2^{-\gamma}$ as $\eta \downarrow 0$.
\end{lemma}

\subsection{Basic harmonic analysis}

Let $\varphi(\xi)$ be a radial bump function supported in the ball
$\{\,\xi \in \mathbb R^d : |\xi| \leq \frac{11}{10}\,\}$ and equal
to 1 on the ball $\{\,\xi \in \mathbb R^d: |\xi |\leq 1 \,\}$. For
each number $N > 0$, we define the Fourier multipliers
\[\aligned
& \widehat{P_{\leq N}f}(\xi)  := \varphi(\xi/N)\hat f(\xi),\\
& \widehat {P_{> N}f}(\xi)  := (1- \varphi(\xi/N))\hat f(\xi), \\
& \widehat{P_N f}(\xi) := \psi(\xi/N)\hat f(\xi)= (\varphi(\xi/N)
- \varphi(2\xi/N))\hat f(\xi)
\endaligned\]
and similarly $P_{<N}$ and $P_{\geq N}$. We also define
\[P_{M < \cdot \leq N} := P_{\leq N} - P_{\leq M} = \sum_{M < N' \leq N}P_{N'}\]
for $M < N$. We will use these multipliers when $M$ and $N$ are
dyadic numbers; in particular, all summations over $N$ or $M$ are
understood to be over dyadic numbers. Nevertheless, it will
occasionally be convenient to allow $M$ and $N$ to not be the
power of 2. Note that, $P_N$ is not truly a projection; to get
around this, define
\[\tilde P_N := P_{N/2} + P_N + P_{2N}.\]
These obey $\tilde P_NP_N = P_N\tilde P_N = P_N $.

The Littlewood-Paley operators commute with the propagator
$e^{it\Delta}$, as well as with differential operators such as
$i\partial_t + \Delta$. We will use basic properties of these
operators many many times. First, we introduce
\begin{lemma}[Bernstein]
For $1 \leq p \leq q \leq \infty$,
\[\aligned
& \big\||\nabla|^{\pm s}P_Nf\big\|_{L_x^q(\mathbb R^d)}  \thicksim
N^{\pm
s}\|P_Nf\|_{L_x^p(\mathbb R^d)}, \\
& \|P_{\leq N}f\|_{L_x^q(\mathbb R^d)}  \lesssim N^{\frac{d}{p}-
\frac{d}{q}}\|P_{\leq N}f\|_{L_x^p(\mathbb R^d)},\\
& \|P_Nf\|_{L_x^q(\mathbb R^d)}  \lesssim N^{\frac{d}{p}-
\frac{d}{q}}\|P_Nf\|_{L_x^p(\mathbb R^d)}.
\endaligned\]
\end{lemma}
We also need the following fractional Leibniz rule, \cite{c31}.
\begin{lemma}[Fractional Leibniz rule]
Let $\alpha \in (0, \, 1),\, \alpha_1,\, \alpha_2 \in [0, \alpha]$
with $\alpha = \alpha_1 + \alpha_2$. Let $1 < p,\, p_1,\,p_2,\, q,
\, q_1, \, q_2 < \infty$ be such that $\frac{1}{p} = \frac{1}{p_1}
+ \frac{1}{p_2}$, $\frac{1}{q} = \frac{1}{q_1} + \frac{1}{q_2}$.
Then
\[\big\| D^\alpha(f g) - g D^\alpha f - f D^\alpha g \big\|_{L_t^qL_x^p}
\lesssim \big\|D^{\alpha_1}
f\big\|_{L_t^{q_1}L_x^{p_1}}\|D^{\alpha_2} g
\|_{L_t^{q_2}L_x^{p_2}}.\] If $\alpha_1 = 0$, $q_1 = \infty $ is
allowed.
\end{lemma}

\subsection{Strichartz's estimates}

Let $e^{it\Delta}$ be  the free Schr\"{o}dinger evolution. From
the explicit formula
\[e^{it\Delta}f(x) = \frac{1}{(4\pi it)^{d/2}} \int_{\mathbb R^d} e^{i|x-y|^2/4t}f(y) \,\mathrm{d} y,\]
we deduce the standard dispersive inequality
\[\|e^{it\Delta}f\|_{L_x^\infty(\mathbb R^d)} \lesssim \frac{1}{|t|^{d/2}}\|f\|_{L_x^1(\mathbb R^d)}\]
for all $t \neq 0$.

 Finer bounds on (frequency localized) linear propagator can be
 derived using  stationary phase:
\begin{lemma}[Kernel estimates, \cite{c13}]
For any $m \geq 0$, the kernel of the linear propagator obeys the
following estimates:
\[|(P_N e^{it\Delta})(x,y)| \lesssim_m \begin{cases}
|t|^{-d/2}, & |x-y| \thicksim N|t| \\
\dfrac{N^d}{|Nt|^m\langle N|x-y|\rangle^m}, & \textrm{otherwise}
\end{cases}\]
for $|t| \geq N^{-2} $ and
\[|(P_Ne^{it\Delta})(x,y)| \lesssim_m N^d \langle N|x-y| \rangle^{-m}\]
for $|t| \leq N^{-2}$.
\end{lemma}

The standard Strichartz's estimate reads:
\begin{lemma}[Strichartz] Let $k \geq 0$, $d \geq 3$. Let $I$ be a
compact time interval, $t_0 \in I$. Then the function $u$ defined
by \begin{equation}\label{e022} u(t) := e^{i(t-t_0)\Delta}u(t_0) -
i \int_{t_0}^t e^{i(t-t')\Delta}f(t') \,\mathrm{d}
t'\end{equation} obeys
\[\|u\|_{\dot S^k(I)} \lesssim \|u(t_0)\|_{\dot H^k_x} + \|f\|_{\dot N^k(I)}\]
for any $t_0 \in I$, where $\dot S^k(I)$ is the Strichartz norm,
and $\dot N^k(I)$ is its dual norm.
\end{lemma}

{\it Proof}. See, for example, \cite{c7}, \cite{c9}. For a
textbook treatment, see \cite{c29}.

 We also
need the following weighted Strichartz's inequality. It is very
useful in regions of space far from the origin.

\begin{lemma}[Weighted Strichartz, \cite{c15}] Let $I$ be an interval, $t_0 \in
I$, $u_0 \in L_x^2(\mathbb R^d)$, $f \in L_{t,x}^{2(d+2)/(d+4)}(I
\times \mathbb R^d)$ be radially symmetric. Then the function $u$
defined by $(\ref{e022})$ obeys the estimate
\[\big\||x|^{\frac{2(d-1)}{q}}u\big\|_{L_t^qL_x^{\frac{2q}{q-4}}(I \times \mathbb R^d)}
\lesssim \|u_0\|_{L_x^2(\mathbb R^d)} +
\|f\|_{L_t^2L_x^{2d/(d+2)}(I \times \mathbb R^d)}\] for all $4
\leq q \leq \infty$.
\end{lemma}

\subsection{In/out decomposition}

For a radially symmetric function $f$, we define the projection
onto outgoing spherical waves by
\[[P^+f](r) = \frac{1}{2}\int_0^\infty r^{\frac{2-d}{2}}H_{\frac{d-2}{2}}^{(1)}(kr) \hat f(k) k^{\frac{d}{2}} \,\mathrm{d} k\]
and the projection onto incoming spherical waves by
\[[P^-f](r) = \frac{1}{2}\int_0^\infty r^{\frac{2-d}{2}}H_{\frac{d-2}{2}}^{(2)}(kr) \hat f(k) k^{\frac{d}{2}} \,\mathrm{d} k\]
where $H_{\frac{d-2}{2}}^{(1)}$ denotes the Hankle function of the
first kind with order $\frac{d-2}{2}$ and
$H_{\frac{d-2}{2}}^{(2)}$ denotes the Hankle function of the
second kind with the same order. We write $P_N^{\pm}$ for the
product $P^{\pm}P_N$, then we have
\begin{lemma}[Kernel estimates, \cite{c15}]
For $|x|\gtrsim N^{-1}$ and $|t| \gtrsim N^{-2}$, the integral
kernel obeys
\[\big|[P_N^{\pm}e^{\mp it\Delta}](x,y)\big| \lesssim \begin{cases}
(|x||y|)^{-\frac{d-1}{2}}|t|^{-\frac{1}{2}}, & |y|-|x| \thicksim
N|t|\\
\dfrac{N^d}{(N|x|)^{\frac{d-1}{2}}\langle N|y|
\rangle^{\frac{d-1}{2}}}\langle N^2t+ N|x|-N|y|\rangle^{-m}, &
\textrm{otherwise} \end{cases}\]
for any $m \geq 0$. For
$|x| \gtrsim N^{-1}$ and $|t| \lesssim N^{-2}$, the integral
kernel obeys
\[\big|[P_N^{\pm}e^{\mp it\Delta}](x,y)\big| \lesssim \frac{N^d}{(N|x|)^{\frac{d-1}{2}}\langle N|y|
\rangle^{\frac{d-1}{2}}}\langle  N|x|-N|y|\rangle^{-m}\] for any
$m \geq 0$.
\end{lemma}
\begin{lemma}[Properties of $P^{\pm}$, \cite{c15}]
We have: \begin{itemize} \item $P^+ + P^-$ acts as the identity on
$L_{rad}^2(\mathbb R^d)$. \item Fix $N > 0$, for any radially
symmetric function $f \in L_x^2(\mathbb R^d)$,\[\|P^{\pm}P_{\geq
N}f\|_{L_x^2(|x| \geq \frac{1}{100N})} \lesssim
\|f\|_{L_x^2(\mathbb R^d)},\]
\end{itemize}
with an $N$-independent constant.
\end{lemma}

\subsection{Concentration compactness}
In this subsection we record the linear profile decomposition
statement due to  Shao
 \cite{c22}. We first recall the symmetries of the solutions to equation $(1.1)$
which fix the initial surface $t =0$.

\begin{definition}[Symmetry group]
For any phase $\theta \in \mathbb R/2\pi\mathbb Z$, position $x_0
\in  \mathbb R^5$, and scaling parameter $\lambda >0$, we define
the unitary transformation $g_{\theta, x_0, \lambda}: \dot
H^{1/2}_x(\mathbb R^5) \mapsto \dot H^{1/2}_x(\mathbb R^5)$ by
\[[g_{\theta,x_0,\lambda}f](x):= \lambda^{-2}e^{i\theta}f(\lambda^{-1}(x-x_0)).\]
Let $G$ denotes the collection of such transformations. For a
function $u : I \times \mathbb R^5 \mapsto \mathbb C$, define
$T_{g_{\theta,x_0,\lambda}}u : \lambda^2I \times \mathbb R^5
\mapsto \mathbb C$ by
\[[T_{g_{\theta,x_0,\lambda}}u](t,x) := \lambda^{-2}e^{i\theta}u(\lambda^{-2}t, \lambda^{-1}(x-x_0))\]
where $\lambda^2I := \{\, \lambda^2t: \, t \in I \, \}$.

Let $G_{rad} \subset G$ denotes the collection of transformations
in $G$ which preserves radial symmetry, or more precisely
\[G_{rad} := \{\,g_{\theta,0,\lambda}: \theta \in \mathbb R/2\pi\mathbb Z, \, \lambda > 0\,\}.\]
\end{definition}
\begin{remark}
$u$ is a maximal life-span solution to $(1.1)$ if and only if
$T_gu$ is a maximal life-span solution to $(1.1)$. Moreover,
$$\|T_gu\|_{\dot H^{1/2}_x(\mathbb R^5)} = \|u\|_{\dot
H^{1/2}_x(\mathbb R^5)}, \quad \|T_gu\|_{S(\lambda^2I)}
=\|u\|_{S(I)}, \quad \textrm{for all} \; \; g \in G.$$
\end{remark}

We are now ready to state the linear profile decomposition.
\begin{lemma}[Linear profiles, \cite{c22}]\label{l29} Let $\{u_n\}_{n\geq 1}$ be a
bounded sequence of functions in $\dot H^{1/2}_x(\mathbb R^5)$.
Then after passing to a subsequence if necessary, there exist a
sequence of functions $\{\phi^j\}_{j\geq 1} \subset \dot
H^{1/2}_x(\mathbb R^5)$, group elements $g_n^j \in G$, and times
$t_n^j \in \mathbb R$ such that we have the decomposition
\begin{equation}\label{e21}
u_n = \sum_{j =1}^J g_n^je^{it_n^j\Delta}\phi^j + \omega_n^J
\end{equation}
for all $J \geq 1$; $\omega_n^J \in \dot H^{1/2}_x(\mathbb R^5)$
obeying
\begin{equation}\label{e22}
\lim_{J \to \infty} \limsup_{n \to
\infty}\|e^{it\Delta}\omega_n^J\|_{L_t^3L_x^{15/4}(\mathbb R
\times \mathbb R^5)} = 0.
\end{equation}
Moreover, for any $j' \neq j$, we have the following orthogonal
property
\begin{equation}\label{e23}
\lim_{n \to \infty} \left( \frac{\lambda_n^j}{\lambda_n^{j'}} +
\frac{\lambda_n^{j'}}{\lambda_n^j} + \frac{|x_n^j -
x_n^{j'}|}{\lambda_n^j} + \frac{|t_n^j -
t_n^{j'}|}{(\lambda_n^j)^2}\right) = 0.
\end{equation}
For any $J \geq 1$
\begin{equation}\label{e24}
\lim_{n \to \infty} \Big[\big\||\nabla|^{\frac{1}{2}}u_n\big\|_2^2
- \sum_{j =1 }^J \big\||\nabla|^{\frac{1}{2}}\phi^j\big\|_2^2 -
\big\||\nabla|^{\frac{1}{2}}\omega_n^J\big\|_2^2\Big] = 0.
\end{equation}
When $\{u_n\}$ is assumed to be radially symmetric, one can choose
$\phi^j, \omega_n^J$ to be radially symmetric and $g_n^j \in
G_{rad}$.
\end{lemma}

The error term also satisfies the following lemma
\begin{lemma}\label{l210}
For all $J \geq 1, \, 1 \leq j \leq J$, the sequence
$e^{-it_n^j\Delta}[(g_n^j)^{-1}\omega_n^J]$ converges weakly to
zero in $\dot H^{1/2}_x(\mathbb R^5)$ as $n \to \infty$.
\end{lemma}
{\it Proof}. The proof is an analogue to  that in \cite{c14},
\cite{c10}.

We end this section with a perturbation theorem
\begin{theorem}
[Long time perturbation theory] Let $I \subset \mathbb R$ be a
compact time interval and let $t_0 \in I$. Let $\tilde u : I
\times \mathbb R^5 \mapsto \mathbb C$ be  a near-solution to
$(1.1)$ in the sense that \[i\partial_t \tilde u + \Delta \tilde u
= F(\tilde u) + e \] for some function $e$.  Suppose $\tilde u$
satisfies
\[\sup_{t \in I}\|\tilde u\|_{\dot H^{1/2}_x(\mathbb R^5)} \leq A,
\quad \|\tilde u\|_{S(I)} \leq M, \quad \|\tilde u\|_{X(I)} <
+\infty, \] for some constant $M,\, A
> 0$. Assume also that
\[\aligned
\|u_0 - \tilde u(t_0)\|_{\dot H^{1/2}_x(\mathbb R^5)} \leq A', \\
\big\||\nabla|^{1/2}e\big\|_{L_t^1L_x^2(I \times \mathbb R^5)}
\leq \varepsilon, \\
\big\|e^{i(t-t_0)\Delta}(u_0 - \tilde u(t_0))\big\|_{S(I)} \leq
\varepsilon.\\
\endaligned\]
 Then, there exists  a solution $u : I \times
\mathbb R^5$ to $(1.1)$ with $u(t_0)=u_0$ such that \[ \sup_{t \in
I}\|u - \tilde u(t)\|_{\dot H^{1/2}_x(\mathbb R^5)} + \|u - \tilde
u\|_{S(I)}  + \|u - \tilde u\|_{X(I)} \leq \varepsilon.\]
\end{theorem}

\section{ Sharp constant for a Hardy-Littlewood-Sobolev type inequality}
\setcounter{equation}{0}

In this section we find the best constant to the following
Hardy-Littlewood-Sobolev type inequality
\begin{equation}\label{a1}
\iint_{\mathbb R^5 \times \mathbb
R^5}\frac{|u(x)|^2|u(y)|^2}{|x-y|^3} \,\mathrm{d} x \,\mathrm{d} y
\leq C_5\big\||\nabla|^{\frac{1}{2}}u\big\|_2^2\big\|\nabla
u\big\|_2^2,
\end{equation}
and obtain a sufficient condition for global existence of equation
$(1.1)$ with initial data in $\dot H^{1}_x(\mathbb R^5) \cap \dot
H^{1/2}_x(\mathbb R^5)$. We find that the best constant $ C_5 =
2\big\||\nabla|^{\frac{1}{2}}Q\big\|_2^{-2}$, where $Q$ is the
solution to $(\ref{e14})$. The approach is essentially from
\cite{c28}.

Consider the Weinstein functional
\[J(u) = \dfrac{\big\||\nabla|^{\frac{1}{2}}u\big\|_2^2 \|\nabla u\|_2^2}{\int_{\mathbb R^5} (|\cdot|^{-3} \ast |u|^2)|u|^2 \,
\mathrm{d}x}\, , \qquad \forall u \in  \dot H^{1}_x(\mathbb R^5)
\cap \dot H^{1/2}_x(\mathbb R^5). \]

First observe that if we set  $u_{a, b} = au(bx)$, then
\[J(u_{a,b}) = J(u), \qquad \big\||\nabla|^{\frac{1}{2}}u_{a,b}\big\|_2^2 = a^2b^{-4}\big\||\nabla|^{\frac{1}{2}}u\big\|_2^2,
\qquad  \|\nabla u_{a,b} \|_2^2 = a^2b^{-3}\|\nabla u\|_2^2.\]
\begin{theorem}\label{t1}
\[C_5^{-1} = \inf_{u \in \dot H^1_x(\mathbb R^5) \cap \dot H^{1/2}_x(\mathbb R^5) \setminus \{0\}}J(u)\]
can be obtained at some $Q \in \dot H^1_x(\mathbb R^5) \cap \dot
H^{1/2}_x(\mathbb R^5)$. In addition, $C_5 =
2\big\||\nabla|^{\frac{1}{2}}Q\|_2^{-2}$.
\end{theorem}

Before proving the theorem, we present some compactness tools.

\begin{lemma}[Radial Lemma]
Let $d \geq 3$, $u \in \dot H^1_{\rm rad}(\mathbb R^d) \cap \dot
H^{1/2}_{\rm rad}(\mathbb R^d)$ be a radially symmetric function.
Then
\begin{equation}\label{a2}
\sup_{x \in \mathbb R^d}|x|^{\frac{2d-3}{4}}|u(x)| \lesssim
\big\||\nabla|^{\frac{1}{2}}u\big\|_2^{\frac{1}{2}}\|\nabla
u\|_2^{\frac{1}{2}}.
\end{equation}
\end{lemma}
{\it Proof.} Suppose first $u \in C^\infty_c(\mathbb R^d)$. We
have
\[
\aligned r^{\frac{2d-3}{2}}u(r)^2 & = -\int_r^\infty \frac{\mathrm
d}{\mathrm{ d} s}\big(s^{\frac{2d-3}{2}}u(s)^2\big)\mathrm{d}s\\
& \leq -2 \int_r^\infty s^{\frac{2d-3}{2}}u(s)u'(s)\mathrm{d}s \\
& \lesssim \big\||x|^{-\frac{1}{2}}u\|_2\big\|\nabla u\|_2,
\endaligned\]
$(\ref{a2})$ follows from Hardy's inequality. The general case
then follows by the density argument.

\begin{lemma}[Compactness Lemma]
\[\dot H^1_{\rm rad}(\mathbb R^d)\cap \dot H^{1/2}_{\rm rad}(\mathbb R^d)
\hookrightarrow L^p(\mathbb R^d) \quad \textrm{for all} \quad
\frac{2d}{d-1} < p < \frac{2d}{d-2}.\]
\end{lemma}
{\it Proof.} Let $\{u_k\}$ be a bounded sequence in $\dot H^1_{\rm
rad}\cap \dot H^{1/2}_{\rm rad}$, then by the weak compactness
principle, there exists $u \in \dot H^1_{\rm rad}\cap \dot
H^{1/2}_{\rm rad}$ such that $u_k \rightharpoonup u$ weakly in
$\dot H^1_{\rm rad}\cap \dot H^{1/2}_{\rm rad}$.

For $\varepsilon > 0$, let $R > 0$ to be chosen later. Given $p$
as in the statement, we have
\[\aligned
\|u_k - u\|_{L^p(\mathbb R^d)} & \leq \|u_k - u\|_{L^p(B_R)} +
\|u_k - u\|_{L^p(\{\,x \,: \,|x| > R\,\})}\\
& \leq \|u_k - u\|_{L^p(B_R)} + \|u_k - u\|_{L^\infty(\{\,x \,:
\,|x|
> R \,\})}^{\frac{p(d-1)-2d}{(d-1)p}}\|u_k - u\|_{L^{\frac{2d}{d-1}}(\mathbb
R^d)}^{\frac{2d}{(d-1)p}}.
\endaligned\]
By  Lemma 3.1, we first choose $R$ large enough so that \[\|u_k -
u\|_{L^\infty(\{\,x\, :\,|x|
> R \,\})}^{\frac{p(d-1)-2d}{(d-1)p}}\|u_k - u\|_{L^{\frac{2d}{d-1}}(\mathbb
R^d)}^{\frac{2d}{(d-1)p}} \leq \frac{\varepsilon}{2}.\] On the
other hand, it follows from Rellich's compactness lemma that
\[\|u_k - u\|_{L^p(B_R)} \leq \frac{\varepsilon}{2}\]
for large $k$ and so $\|u_k - u\|_{L^p(\mathbb R^d)} \leq
\varepsilon$. This proves the lemma.

\vspace{8pt}

 {\it Proof of Theorem $\ref{t1}$ }. Since $J(u) \geq 0$, we may find a
minimizing sequence $\{u_k\} \subset \dot H^1 \cap \dot H^{1/2}$
such that
\[C_5^{-1} = \inf J(u) = \lim_{k \to \infty}J(u_k) .\]
By symmetric rearrangement technique, we may assume $u_k > 0$ and
is radially symmetric for all $k$.

Set $a_k = \big\||\nabla|^{\frac{1}{2}}u_k\big\|_2^3/{\|\nabla
u_k\|_2^4}$, $b_k =
\big\||\nabla|^{\frac{1}{2}}u_k\big\|_2^2/{\|\nabla u_k\|_2^2}$,
and $Q_k = a_k u(b_k x)$. Then $Q_k \geq 0$, is radially
symmetric. Moreover, we have
\[\big\||\nabla|^{\frac{1}{2}}Q_k\big\|_2 = \|\nabla Q_k\|_2 =1, \quad \lim_{k \to \infty} J(Q_k) = C_5^{-1}.\]
Since $\{Q_k\} \subset \dot H^1_{\rm rad} \cap \dot H^{1/2}_{\rm
rad}$ is uniformly bounded, up to a subsequence,  $Q_k
\rightharpoonup Q^*$ in $\dot H^1_{\rm rad} \cap \dot H^{1/2}_{\rm
rad}$, and $\big\||\nabla|^{\frac{1}{2}}Q^*\big\|_2 \leq 1$,
$\|\nabla Q^*\|_2 \leq 1$. From Lemma $3.2$, $Q_k \to Q^*$ in
$L^p(\mathbb R^5)$ for $\frac{5}{2} < p < \frac{10}{3}$.
Furthermore, we have
\[\iint_{\mathbb R^5 \times \mathbb R^5} \frac{|Q_k(x)|^2|Q_k(y)|^2}{|x - y|^3} \, \mathrm{d}x \mathrm{d}y
\longrightarrow \iint_{\mathbb R^5 \times \mathbb R^5}
\frac{|Q^*(x)|^2|Q^*(y)|^2}{|x - y|^3} \, \mathrm{d}x \mathrm{d}y
\quad \textrm{as} \, k \to \infty.\] This is easily checked by a
direct computation using the Hardy-Littlewood-Sobolev inequality.

Thus
\[C_5^{-1} \leq  J(Q^*) \leq \frac{1}{\int_{\mathbb R^5}(|\cdot|^{-3} \ast |Q^*|^2)|Q^*|^2\, \mathrm{d}x} =
\lim_{k \to \infty}J(Q_k) = C_5^{-1}.\] This implies that
$\big\||\nabla|^{\frac{1}{2}}Q^*\big\|_2^2\|\nabla Q^*\|_2^2 =1$,
which further gives $\big\||\nabla|^{\frac{1}{2}}Q^*\big\|_2 =
\|\nabla Q^*\|_2 =1$.

Since $Q^*$ is a minimizer, it satisfies the Euler-Lagrangian
equation
\[\frac{\mathrm{d}}{\mathrm{d}\varepsilon}\Big|_{\varepsilon =0}J(Q^* + \varepsilon \phi) = 0 \quad \textrm{for all}\, \phi \in
C_0^\infty(\mathbb R^5).\]
 Taking into account the fact that $\big\||\nabla|^{\frac{1}{2}}Q^*\big\|_2 = \|\nabla Q^*\|_2
=1$, we have
\[-\Delta Q^* + (-\Delta)^{1/2}Q^* - 2C_5^{-1}(|\cdot|^{-3} \ast |Q^*|^2)Q^* = 0.\]
Let $Q^* = \sqrt{C_5/2} Q$, then $Q$ solves $(\ref{e14})$.

By the fact that $ \big\||\nabla |^{\frac{1}{2}}Q^*\big\|_2 =1$,
it yields $C_5 =2 \big\||\nabla|^{\frac{1}{2}}Q\big\|_2^{-2}$.
$\square$

\begin{proposition}\label{t2}
Let $u_0 \in \dot H^1_x(\mathbb R^5) \cap \dot H^{1/2}_x(\mathbb
R^5)$. Suppose $\sup_t \||\nabla|^{\frac{1}{2}}u\|_2 <
\||\nabla|^{\frac{1}{2}} Q\|_2$, then the solution to $(1.1)$ is
global.
\end{proposition}
{\it Proof.} It is a consequence of the energy conservation
\[E(u(t)) = \frac{1}{2}\int_{\mathbb R^5}|\nabla u|^2 \, \mathrm{d}x - \frac{1}{4}\iint_{\mathbb R^5 \times \mathbb R^5}
\frac{|u(x)|^2|u(y)|^2}{|x-y|^3} \, \mathrm{d}x \mathrm{d}y,\] and
$(\ref{a1})$.

\section{Reduction to almost periodic solution}
\setcounter{equation}{0}

\vspace{12pt}

In this section we will prove Theorem 1.3. The main step toward
this end is to prove a Palais-Smale condition modulo scaling.

For any $A > 0$, define
\[ L(A) = \sup\left\{\|u\|_{S(I)}:\,\, u: I \times \mathbb R^5 \mapsto \mathbb C \, \textrm{such that }\,
 \sup_{t \in I} \|u\|_{\dot H^{1/2}_x} \leq A\right\}.\]
Here, the supremum is taken over all solutions $u : I \times
\mathbb R^5 \mapsto \mathbb C$ to (1.1) satisfying $\sup_{t \in I}
\|u\|_{\dot H^{1/2}_x} \leq A$. Note that $L(A)$ is non-decreasing
and left-continuous. On the other hand, from Theorem 1.1,
\[ L(A) \lesssim A  \quad \textrm{for } \quad A \leq \delta_0,\]
where $\delta_0$ is the threshold from the small data global
well-posedness theory. Theorem 1.2 states that for each $A <
\frac{\sqrt{6}}{3} \|Q\|_{\dot H^{1/2}}$, $L(A) < \infty$.
Therefore, if Theorem 1.2 failed,  there exists $\delta_0 < A_c <
\frac{\sqrt{6}}{3}\|Q\|_{\dot H^{1/2}}$ such that $L(A) < +\infty$
for $A < A_c$, $L(A) = +\infty$ for $A \geq A_c$.

{\it Convention:} In this section and the rest sections, we write
$|x|^{-3} \ast$ as $|\nabla|^{-2}$ since they are equivalent up to
a constant. Moreover, we ignore the distinction between a function
and its conjugation as they make no difference in our discussion.

\subsection{Palais-Smale condition modulo scaling}
\begin{proposition}\label{p41}
Let $u_n: I_n \times \mathbb R^5 \mapsto \mathbb C$ be a sequence
of solutions to $(1.1)$ such that
\begin{equation}\label{e31}
\limsup_{n \to \infty} \sup_{t \in I_n} \|u_n(t)\|_{\dot
H^{1/2}_x} = A_c.
\end{equation}
Let $t_n \in I_n$ be a time sequence such that
\[\lim_{n \to \infty}\|u_n\|_{S(-\infty,\; t_n)} = \lim_{n \to \infty}\|u_n\|_{S(t_n,\; \infty)} = \infty.\]
Then there exists a  subsequence  of  $u_n(t_n)$,  which converges
in $\dot H^{1/2}_x(\mathbb R^5)$ modulo scaling.
\end{proposition}
The proof of this Proposition is achieved through several steps.

{\it Proof}. By time-translation invariant of (1.1), we may set
$t_n = 0$ for all $n \geq 1$. Then
\begin{equation}\label{e32}
\lim_{n \to \infty}\|u_n\|_{S(-\infty, \; 0)} = \lim_{n \to
\infty}\|u_n\|_{S(0, \; \infty)} = \infty.
\end{equation}
Now applying Lemma $\ref{l29}$ to the sequence $u_n(0)$, and up to
a subsequence, we obtain a decomposition
\[u_n(0) = \sum_{j= 1}^J g_n^je^{it_n^j\Delta}\phi^j  + \omega_n^J\]
for any $J \geq 1$, $n \geq 1$.

 By passing to a further subsequence, we may assume   $t_n^j$ converges to some $t^j \in
[-\infty,+\infty]$ for each $j$. If $t^j$ is finite, then
replacing $\phi^j$ by $e^{it^j\Delta}\phi^j$, we may set $t^j =0$.
Adding $e^{it_n^j\Delta}\phi^j - \phi^j$ to the error term
$\omega_n^J$, we may assume $t_n^j \equiv 0$ . Thus, we only need
to deal with $t_n^j \equiv 0$ and $t_n^j \to \pm \infty$.

For each $\phi^j$ and $t_n^j$, define  nonlinear profile $v^j :
I^j \times \mathbb R^5 \mapsto \mathbb C$  as follows:
\begin{itemize}
\item If $t_n^j \equiv 0$, then $v^j$ is the maximal life-span
solution to (1.1) with initial data $v^j(0) = \phi^j$.  \item If
$t_n^j \to \infty$, then $v^j$ is the maximal life-span solution
to (1.1) that scatters forward to $e^{it\Delta}\phi^j$. \item If
$t_n^j \to -\infty$, then $v^j$ is the maximal life-span solution
to (1.1) that scatters backward to $e^{it\Delta}\phi^{j}$.
\end{itemize}
For each $j$, $n \geq 1$, define $v_n^j : I_n^j \times \mathbb R^5
\mapsto \mathbb C$ by
\[v_n^j(t) := T_{g_n^j}[v^j(\cdot + t_n^j)](t),\]
where $I_n^j :=\{\,t \in \mathbb R:  \, (\lambda_n^j)^{-2}t +
t_n^j \in I^j\,\}$. Then for each $j$, $v_n^j$ is also a maximal
life-span solution to (1.1) with initial data $v_n^j(0) =
g_n^jv^j(t_n^j)$, and with maximal life-span $I_n^j=(-T^-_{n,j},\;
T^+_{n,j})$, $-\infty \leq -T_{n,j}^- < 0 < T_{n,j}^+ \leq
+\infty$.

With these preliminaries out of the way, we first have

{\bf Step 1:} There exists $J_0 \geq 1$ such that, for all  $j \geq
J_0$, $n$ sufficiently large
\begin{equation}\label{e33}
\sup_{t \in \mathbb R} \|v_n^j(t)\|_{\dot H^{1/2}_x} +
\|v_n^j\|_{S(\mathbb R)} + \|v_n^j\|_{X(\mathbb R)}\lesssim
\|\phi^j\|_{\dot H^{1/2}_x}.
\end{equation}
{\it Proof}. From $(\ref{e24})$, there exists $J_0 \geq 1$ such
that for sufficiently large $n$
$$\|\phi^j\|_{\dot H^{1/2}_x} \leq \delta_0 \quad \textrm{for all} \quad j \geq J_0$$
 where $\delta_0 $
is the threshold from the small data theory. Hence, by Theorem
1.1, $v_n^j$ is global and \[\sup_{t \in \mathbb R}\|v_n^j\|_{\dot
H^{1/2}_x} +\|v_n^j\|_{X(\mathbb R)} + \|v_n^j\|_{S(\mathbb R)}
\lesssim \|\phi^j\|_{\dot H^{1/2}_x}.\] for all  $j \geq J_0$ and
all $n$ sufficiently large.

\vspace{8pt}

{\bf Step 2:} There exists $1 \leq j_0 < J_0$ such that
\[\limsup_{n \to \infty} \|v_n^{j_0}\|_{S(0,\; T_{n,j_0}^+)} = \infty.\]
{\it Proof}. Suppose to the contrary that for all $1 \leq j < J_0$
\begin{equation}\label{e34}
\limsup_{n \to \infty} \|v_n^j\|_{S(0, \; T_{n,j}^+)} \leq M <
\infty
\end{equation}
for some $M > 0$. This implies that $T_{n,j}^+ = \infty$ for all
$1 \leq j < J_0$ and all sufficiently large $n$. Given $\eta > 0$,
divide $(0, \infty)$ into subintervals $I_k$ such that on each
$I_k$, $\|v_n^j\|_{S(I_k)} \leq \eta$. By Strichartz's estimate,
we have for all $1 \leq j < J_0$ and all large $n$ that
\begin{equation}\label{e35}
\|v_n^j\|_{X(0,\infty)} < \infty.
\end{equation}
Indeed, let $\eta > 0$, divide $(0, \infty)$ into subintervals
$I_k = [t_k, t_{k+1}]$ such that on each $I_k$ we have
$\|v_n^j\|_{S(I_k)} \leq \eta$. Note that, there are at most
$\eta^{-1} \times M$ such intervals. Applying the Strichartz
estimate
\begin{eqnarray*}
\|v_n^j\|_{X(I_k)}& \lesssim & \|v^j_n(t_k)\|_{\dot H^{1/2}_x} +
\big\||\nabla|^{\frac{1}{2}}F(v_n^j)\big\|_{L_t^1 L_x^2} \nonumber\\
& \lesssim & A_c + \|v_n^j\|^2_{S(I_k)}\|v_n^j\|_{X(I_k)}.
\end{eqnarray*}
If we choose $\eta > 0$ sufficiently small, then
\[\|v_n^j\|_{X(I_k)} \lesssim A_c.\]
Summing over all $I_k$, we achieve $(\ref{e35})$.

Combining $(\ref{e34})$ with Step 1, and then using $(\ref{e24})$
and $(\ref{e31})$, we have that for all sufficiently large $n$,
\begin{equation}\label{e36}
\sum_{j \geq 1} \sup_{t \in (0,\infty)}\|v_n^j\|_{\dot H^{1/2}_x}
+ \|v_n^j\|_{S(0,\infty)} + \|v_n^j\|_{X(0, \infty)} \lesssim 1 +
A_c.
\end{equation}

Next, we will use perturbation theorem to obtain a bound on
$\|u_n\|_{S(0,\, \infty)}$ for $n$ sufficiently large.

Define an approximation to $u_n$ by
\begin{equation}\label{e00}
u_n^J(t) := \sum_{j =1}^J v_n^j (t) + e^{it\Delta}\omega_n^J.
\end{equation}
Then, by the definition of nonlinear profile
\[\aligned
\limsup_{n \to \infty}\|u_n^J(0) - u_n(0)\|_{\dot H^{1/2}_x} & =
\limsup_{n \to \infty}\Big\|\sum_{j=1}^J g_n^jv^j(t_n^j) -
g_n^je^{it_n^j\Delta}\phi^j\Big\|_{\dot H^{1/2}_x}
\\
& \lesssim \limsup_{n \to \infty} \sum_{j=1}^J \|v^j(t_n^j) -
e^{it_n^j\Delta}\phi^j\|_{\dot H^{1/2}_x} = 0.
\endaligned\]

Note that  $(\ref{e23})$ with a few computations yields that for
all $j \geq 1$
\begin{equation}\label{e37}
\limsup_{n \to \infty }\|v_n^{j'}v_n^j\|_{S(0, \; \infty)} =0
\quad
\end{equation}
for any $ j' \neq j$.(Such an asymptotic orthogonal property was
well developed in \cite{c12}, \cite{c22}, we refer to them for
details.)

Thus, by $(\ref{e22})$,  $(\ref{e36})$ and $(\ref{e37})$ \vspace{-12pt}
\begin{eqnarray}\label{e38}
\lim_{J \to \infty}\limsup_{n \to \infty}\|u_n^J\|_{S(0, \;
\infty)} & \lesssim &\lim\limits_{J \to \infty}\limsup\limits_{n
\to \infty} \Big(\Big\|\sum_{j=1}^J v_n^j\Big\|_{S(0, \; \infty)}
+
\big\|e^{it\Delta}\omega_n^J\big\|_{S(0, \; \infty)}\Big) \nonumber\\
& \lesssim & \lim\limits_{J\to \infty} \limsup\limits_{n\to
\infty} \sum_{j=1}^J \|v_n^j\|_{S(0, \; \infty)} \lesssim 1 + A_c.
\end{eqnarray}

By the same argument as that  to derive $(\ref{e35})$ from
$(\ref{e34})$, we obtain
\[
\lim_{J \to \infty} \limsup_{n \to \infty} \|u_n^J\|_{X(0, \;
\infty)} < \infty.
\]

Now, we have to verify that
\[\lim_{J \to \infty} \limsup_{n \to \infty}
\Big\||\nabla|^{\frac{1}{2}}\big[(i\partial_t +\Delta)u_n^J +
F(u_n^J)\big]\Big\|_{L_t^1L_x^2((0,\infty)\times \mathbb R^5)}
=0.\]

Using the triangle inequality, we need to show on $(0,\infty)
\times \mathbb R^5$ that
\begin{equation}\label{e39}
 \lim_{J \to \infty} \limsup_{n \to \infty} \Big\||\nabla|^{\frac{1}{2}}\Big[\sum_{j=1}^J F(v_n^j)- F(\sum_{j=1}^J
v_n^j)\Big]\Big\|_{L_t^1L_x^2} =0\end{equation} and
\begin{equation}\label{e310}
 \lim_{J \to \infty} \limsup_{n \to \infty} \Big\||\nabla|^{\frac{1}{2}}\big(F(u_n^J - e^{it\Delta}\omega_n^J) -
F(u_n^J)\big)\Big\|_{L_t^1L_x^2} = 0. \end{equation}

We first consider $(\ref{e39})$. By expanding out the nonlinearity
\[
\aligned & \Big||\nabla|^{\frac{1}{2}}\Big[\sum_{j=1}^J F(v_n^j)-
F(\sum_{j=1}^J v_n^j)\Big] \Big|\\
 \leq & \sum_{j_1, j_2,j_3=1}^J
 \Big||\nabla|^{\frac{1}{2}}\big[\big(|\nabla|^{-2}
 (v_n^{j_1}{v_n^{j_2}})\big)v_n^{j_3}\big]\Big|,
\endaligned\]
where at least two of $j_1, j_2, j_3$ are different.

Note that the nonlocal action (i.e. convolution) break up the
spatial orthogonality, whereas time orthogonality will be
preserved. Recalling the radial assumption, we may assume $j_2
\neq j_1$. Thus, using the fractional Leibniz rule, H\"{o}lder's
inequality,the  Hardy-Littlewood-Sobolev inequality, and
$(\ref{e37})$, we obtain on
 $(0,\infty) \times \mathbb R^5$ that
\[
\aligned & \lim_{J \to \infty}\limsup_{n \to \infty}
\Big\||\nabla|^{\frac{1}{2}}\Big[\sum_{j=1}^J F(v_n^j)-
F(\sum_{j=1}^J
v_n^j)\Big]\Big\|_{L_t^1L_x^2} \\
\lesssim_J & \lim_{J\to \infty} \limsup_{n \to
\infty}\sum_{j_1,j_2,j_3=1}^J \Big(
\big\||\nabla|^{\frac{1}{2}}\big(|\nabla|^{-2}
(v_n^{j_1}{v_n^{j_2}})\big)v_n^{j_3}\big\|_{L_t^1L_x^2} +
\big\|\big(|\nabla|^{-2}
(v_n^{j_1}{v_n^{j_2}})\big)|\nabla|^{\frac{1}{2}}v_n^{j_3}\big\|_{L_t^1L_x^2}\Big)\\
\lesssim_J & \lim_{J\to \infty} \limsup_{n \to
\infty}\sum_{j_1,j_2,j_3=1}^J \Big(
\big\||\nabla|^{\frac{1}{2}}\big(|\nabla|^{-2} (v_n^{j
_1}{v_n^{j_2}})\big)\big\|_{L_t^{\frac{3}{2}}L_x^{\frac{30}{7}}}\|v_n^{j_3}\|_{S(0,
\infty)}
 \\
& \hspace{4cm} + \big\||\nabla|^{-2}
(v_n^{j_1}{v_n^{j_2}})\big\|_{L_t^{\frac{3}{2}}L_x^{\frac{15}{2}}}\|v_n^{j_3}\|_{X(0, \infty)} \Big)\\
 \lesssim_{J}& \lim_{J \to \infty}\limsup_{n \to
\infty}\sum_{j_1,j_2,j_3=1}^J
\|v_n^{j_1}{v_n^{j_2}}\|_{L_t^{\frac{3}{2}}L_x^{\frac{15}{8}}} =
0, \endaligned\] where  the last limit is also a consequence of
the orthogonality.

For $(\ref{e310})$, note that on $(0,\infty) \times \mathbb R^5$
\begin{eqnarray*}
&  &\big\||\nabla|^{\frac{1}{2}}(F(u_n^J -
e^{it\Delta}\omega_n^J) - F(u_n^J))\big\|_{L_t^1L_x^2}  \\
& \lesssim & \big\||\nabla|^{\frac{1}{2}}[\big(|\nabla|^{-2} (
u_n^J {e^{it\Delta}\omega_n^J})\big)u_n^J]\big\|_{L_t^1L_x^2} +
\big\||\nabla|^{\frac{1}{2}}[(|\nabla|^{-2}(u_n^J
{e^{it\Delta}\omega_n^J}))e^{it\Delta}\omega_n^J]\big\|_{L_t^1L_x^2}\\
& & + \big\||\nabla|^{\frac{1}{2}}[(|\nabla|^{-2}
|u_n^J|^2)e^{it\Delta}\omega_n^J]\big\|_{L_t^1L_x^2} +
\big\||\nabla|^{\frac{1}{2}}[(|\nabla|^{-2}
|e^{it\Delta}\omega_n^J|^2)e^{it\Delta}\omega_n^J]\big\|_{L_t^1L_x^2}\\
 & & + \big\||\nabla|^{\frac{1}{2}}[(|\nabla|^{-2}
|e^{it\Delta}\omega_n^J|^2)u_n^J]\big\|_{L_t^1L_x^2}.
\end{eqnarray*}
 Using $(\ref{e22})$, H\"{o}lder's  inequality, the Hardy-Littlewood-Sobolev inequality, the above terms  on the right hand side will go to
 zero as $J$, $n$ tend to $\infty$, except
\[\big\||\nabla|^{\frac{1}{2}}[(|\nabla|^{-2}
|u_n^J|^2)e^{it\Delta}\omega_n^J]\big\|_{L_t^1L_x^2((0,\;
\infty)\times \mathbb R^5)}.\] By the fractional Leibniz rule and
the triangle inequality, it suffices to estimate
\[\big\||\nabla|^{\frac{1}{2}}(|\nabla|^{-2}  |u_n^J|^2)e^{it\Delta}\omega_n^J\big\|_{L_t^1L_x^2((0,\;\infty)\times \mathbb R^5)}\]
and
\[\big\|(|\nabla|^{-2} |u_n^J|^2)|\nabla|^{\frac{1}{2}}e^{it\Delta}\omega_n^J\big\|_{L_t^1L_x^2((0,\; \infty)\times \mathbb R^5)}.\]
Using H\"{o}lder's, the Hardy-Littlewood-Sobolev inequality, and
$(\ref{e22})$, the first integral goes to zero when $J$, $n$ go to
infinity. Then, we are reduced to showing that the second integral
has limit zero with $J$, $n$.

Replace $u_n^J$ with its definition formula $(\ref{e00})$ to get
on $(0,\infty) \times \mathbb R^5$
\[
\aligned &\big\|(|\nabla|^{-2}
|u_n^J|^2)|\nabla|^{\frac{1}{2}}e^{it\Delta}\omega_n^J\big\|_{L_t^1L_x^2}\\
\lesssim & \sum_{j=1}^J \big\|(|\nabla|^{-2}
|v_n^j|^2)|\nabla|^{\frac{1}{2}}e^{it\Delta}\omega_n^J\big\|_{L_t^1L_x^2}
+ \sum_{j' \neq j} \big\|(|\nabla|^{-2}(
v_n^j{ v_n^{j'}}))|\nabla|^{\frac{1}{2}}e^{it\Delta}\omega_n^J\big\|_{L_t^1L_x^2}\\
& + \sum_{j=1}^J \big\|(|\nabla|^{-2}( v_n^j
{e^{it\Delta}\omega_n^J}))|\nabla|^{\frac{1}{2}}e^{it\Delta}\omega_n^J\big\|_{L_t^1L_x^2}
:= {\rm I_1 + I_2 + I_3}.
\endaligned\]
By $(\ref{e23})$, ${\rm I}_2$ will go to zero as $J$, $n$ go to
infinity. Using $(\ref{e22})$, ${\rm I}_3$ vanishes as $J$, $n$
tend to infinity. So, We only need to show that ${\rm I}_1$ also
vanishes.

For arbitrary $\eta > 0$, from $(\ref{e38})$, there exists
$J'(\eta) \geq 1 $ such that
\[\sum_{j \geq J'} \|v_n^j\|_{S(0, \; \infty)} \leq \eta.\]
Thus, we are reduced to proving that
\[\lim_{J \to \infty}\limsup_{n\to \infty} \big\|(|\nabla|^{-2}
|v_n^j|^2)|\nabla|^{\frac{1}{2}}e^{it\Delta}\omega_n^J\big\|_{L_t^1L_x^2((0,
\; \infty)\times \mathbb R^5)} =0 \quad \textrm{for all} \quad 1
\leq j \leq J'.\]

Fix $1 \leq j \leq J'$. A change of variables yields
\[ \big\|(|\nabla|^{-2}
|v_n^j|^2)|\nabla|^{\frac{1}{2}}e^{it\Delta}\omega_n^J\big\|_{L_t^1L_x^2}
 =  \big\|\big(|\nabla|^{-2}
|v^j|^2\big)|\nabla|^{\frac{1}{2}}\big[T_{(g_n^j)^{-1}}(e^{it\Delta}\omega_n^J)\big](\cdot
- t_n^j)\big\|_{L_t^1L_x^2}.
\]
Let $\tilde \omega_n^J :=
[T_{(g_n^j)^{-1}}(e^{it\Delta}\omega_n^J)](\cdot - t_n^j)$,
$\mathcal{ I} : v^j \mapsto (|\nabla|^{-2} |v^j|^2)$. Note that
\begin{equation}\label{e311}
\|\tilde \omega_n^J\|_{S(0,\;\infty)} =
\|e^{it\Delta}\omega_n^J\|_{S(0,\infty)}, \quad \|\tilde
\omega_n^J\|_{X(0,\infty)} = \|e^{it\Delta}\omega_n^J\|_{X(0,\;
\infty)}.
\end{equation}
Using H\"{o}lder's inequality, the interpolation theorem, we see
\[
\aligned & \big\|\mathcal {I}(v^j)|\nabla|^{\frac{1}{2}}\tilde
\omega_n^J \big\|_{L_t^1L_x^2} \\
\lesssim & \|\mathcal
{I}(v^j)\|_{L_t^{12/7}L_x^{15}}\big\||\nabla|^{\frac{1}{2}}\tilde
\omega_n^J\big\|_{L_t^{12/5}L_x^{30/13}} \\
\lesssim & \|v^j\|_{L_t^{24/7}L_x^{30/7}}\big\|\tilde
\omega_n^J\big\|_{X(0,
\infty)}^{1/2}\big\||\nabla|^{\frac{1}{2}}\tilde
\omega_n^J\big\|_{L_{t,x}^2}^{1/2}.
\endaligned\]
By density, we may assume $\mathcal I(v_n^j) \in
C_c^\infty(\mathbb R \times \mathbb R^5)$. It thus suffices to
verify
\[\lim_{J \to \infty}\limsup_{n\to \infty}\big\||\nabla|^{\frac{1}{2}}\tilde \omega_n^J\big\|_{L_{t,x}^2(K)} = 0\]
for any compact $K \subset \mathbb R \times \mathbb R^5$. This is
a consequence of  $(\ref{e22})$ and the following lemma:
\begin{lemma}
Let $\phi \in \dot H^{1/2}_x(\mathbb R^5)$. Then
\[\big\||\nabla|^{\frac{1}{2}}e^{it\Delta}\phi\big\|_{L_{t,x}^2(\,[-T,\;T]\times \{\,x : \,|x| \leq R\,\}\,)}^2
\lesssim T^{\frac{1}{6}}R^{\frac{5}{3}}\|e^{it\Delta}\phi\|
_{L_t^3L_x^{15/4}}\big\||\nabla|^{\frac{1}{2}}\phi\big\|_{L_x^2}.\]
\end{lemma}
{\it Proof.} The proof is  analogous  to the one of Lemma 2.5 in
\cite{c14}.

 Now, applying perturbation theorem  with $\tilde u = u_n^J$, $e
= (i\partial_t + \Delta)u_n^J - F(u_n^J)$, and using $(\ref{e38})$,
we obtain
\[\|u_n^J\|_{S(0, \;\infty)} \lesssim 1 + A_c\]
for all sufficiently large $n$. This contradicts $(\ref{e32})$,
which concludes Step 2.

\vspace{12pt}
 Combining Step 1
with Step 2, and rearranging the indices, we may find $1 \leq J_1
\leq J_0$ such that
\[
\aligned & \limsup_{n \to \infty} \|v_n^j\|_{S(0,\; T_{n,j}^+)} =
\infty \quad
\textrm{for} \,\, 1 \leq j \leq J_1, \\
& \limsup_{n \to \infty} \|v_n^j\|_{S(0, \; T_{n,j}^+)} < \infty
\quad \textrm{for} \,\,   j > J_1.
\endaligned
\]

For $m \in \mathbb N$, $n \geq 1$,  define an interval $K_n^m$ of
the form $[0,\tau]$ by
\[\sup_{1 \leq j \leq J_1} \|v_n^j\|_{S(K_n^m)} = m .\]
Then, $v_n^j$ is defined on $K_n^m$ for all $j \geq 1$ and
$\|v_n^j\|_{S(K_n^m)}$ is finite for all $j \geq 1$.

Since $u_n^J$ is a good approximation to $u_n$, using the same
argument as in Step 2, we may obtain

\begin{equation}\label{e313}
\lim_{J \to \infty}\limsup_{n \to \infty}\sup_{t \in K_n^m} \|u_n^J
- u_n\|_{\dot H^{1/2}_x(\mathbb R^5)} =0
\end{equation}
for each $m \geq 1$.

By the definition of  $K_n^m$, we may choose $1 \leq j_0=j_0(m,n)
\leq J_1$ such that
\begin{equation}
\|v_n^{j_0(m,n)}\|_{S(K_n^m)} = m.
\end{equation}
Moreover,  there are infinitely many $m$ satisfying $j_0(m,n) = j_0$
for infinitely many $n$.

By the definition of $A_c$, we have
\begin{equation}\label{e312}
\limsup_{m\to \infty}\limsup_{n\to \infty}\sup_{t \in
K_n^m}\|v_n^{j_0}\|_{\dot H^{1/2}_x(\mathbb R^5)} \geq A_c.
\end{equation}

{ \bf Step 3:} For all $J\geq 1$ and $m \geq 1$
\begin{equation}\label{e314}
\lim_{n \to \infty} \sup_{t \in K_n^m} \Big|\|u_n^J(t)\|^2_{\dot
H^{1/2}_x} - \sum_{j =1}^J \|v_n^j(t)\|^2_{\dot H^{1/2}_x} -
\|\omega_n^J\|^2_{\dot H^{1/2}_x} \Big| = 0.
\end{equation}
{\it Proof}.  Note that for all $J \geq 1$, $m \geq 1$
\[
\aligned \|u_n^J(t)\|_{\dot H^{1/2}_x}^2 & = \big\langle \,
|\nabla|^{\frac{1}{2}}u_n^J(t), |\nabla|^{\frac{1}{2}}u_n^J(t)\,
\big\rangle\\
& = \sum_{j =1}^J \big\||\nabla|^{\frac{1}{2}}v_n^j\big\|_{\dot
H^{1/2}_x}^2 + \|\omega_n^J\|^2_{\dot H^{1/2}_x} + \sum_{j' \neq
j} \big\langle\, |\nabla|^{\frac{1}{2}}v_n^j(t),
|\nabla|^{\frac{1}{2}}v_n^{j'}(t)\,\big\rangle \\
&  \hspace{24pt}+ \sum_{j=1}^J \Big(\big\langle\,
|\nabla|^{\frac{1}{2}}e^{it\Delta}\omega_n^J,
|\nabla|^{\frac{1}{2}}v_n^j(t)\,\big\rangle + \big\langle\,
|\nabla|^{\frac{1}{2}}v_n^j(t),|\nabla|^{\frac{1}{2}}e^{it\Delta}\omega_n^J
\,\big\rangle\Big).
\endaligned\]
Thus, to establish $(\ref{e314})$, it suffices to show that for all
$t_n \in K_n^m$,
\begin{equation}\label{e315}
\lim_{n \to \infty}\big\langle \,|\nabla|^{\frac{1}{2}}v_n^j(t_n)
\, , \, |\nabla|^{\frac{1}{2}}v_n^{j'}(t_n)\, \big\rangle = 0
\end{equation}
and
\begin{equation}\label{e316}
\lim_{n \to \infty}\big\langle\,
|\nabla|^{\frac{1}{2}}e^{it_n\Delta}\omega_n^J \, , \,
|\nabla|^{\frac{1}{2}}v_n^{j}(t_n) \,\big\rangle = 0
\end{equation}
for all $1 \leq j, \, j' \leq J$, $j \neq j'$.

We only deal with $(\ref{e316})$, as $(\ref{e315})$ can be done in
the same manner, using $(\ref{e23})$.

Do a change of variables, the formula in $(\ref{e316})$ becomes
\begin{equation}\label{e317}
\big\langle\,
|\nabla|^{\frac{1}{2}}e^{it_n(\lambda_n^j)^{-2}\Delta}[(g_n^j)^{-1}\omega_n^J]\,
, \, |\nabla|^{\frac{1}{2}}v^j(t_n^j + t_n(\lambda_n^j)^{-2})
\,\big\rangle.
\end{equation}
Since $t_n \in K_n^m \subset [0,T_{n,j}^+)$ for all $1 \leq j \leq
J_1$ and $v_j$ has maximal-life span $I^j = \mathbb R$ for $j >
J_1$, we have $t_n(\lambda_n^j)^{-2} + t_n^j \in I^j$ for all $j
\geq 1$. By passing to a subsequence in $n$, we may assume
$t_n(\lambda_n^j)^{-2} + t_n^j \to \tau^j $.

If $\tau^j$ is finite, then by the continuity of the flow,
$v^j(t_n(\lambda_n^j)^{-2}+t_n^j) \to v^j(\tau^j)$ in $\dot
H^{1/2}_x$.

From $(\ref{e24})$, we have
\begin{eqnarray*}
& & \lim_{n \to \infty}
\big\|e^{it_n(\lambda_n^j)^{-2}\Delta}[(g_n^j)^{-1}\omega_n^J]\big\|_{\dot
H^{1/2}_x(\mathbb R^5)}  = \lim_{n \to \infty}\|\omega_n^J\|_{\dot
H^{1/2}_x} \lesssim A_c.
\end{eqnarray*}
Combining this with $(\ref{e317})$, and using Lemma $\ref{l210}$,
we obtain
\[
\aligned &\lim_{n \to \infty}\big\langle\,
|\nabla|^{\frac{1}{2}}e^{it_n\Delta}\omega_n^J\, , \,
|\nabla|^{\frac{1}{2}}v_n^j(t_n^j)\,\big\rangle\\
= & \lim_{n \to \infty}\big\langle\,
|\nabla|^{\frac{1}{2}}e^{it_n(\lambda_n^j)^{-2}\Delta}[(g_n^j)^{-1}\omega_n^J]\,
, \,
|\nabla|^{\frac{1}{2}}v^j(\tau^j)\,\big\rangle\\
=& \lim_{n \to \infty}\big\langle\,
|\nabla|^{\frac{1}{2}}e^{-it_n^j\Delta}[(g_n^j)^{-1}\omega_n^J]\,
, \,
|\nabla|^{\frac{1}{2}}e^{-i\tau^j\Delta}v^j(\tau^j)\,\big\rangle\\
= &  \, \, 0,
\endaligned\]
which concludes $(\ref{e314})$.

If $\tau^j = +\infty$, then since $t_n(\lambda_n^j)^{-2} \geq 0$,
we must have $\sup I^j = \infty$ and $v^j$ scatters forward in
time. Therefore, there exists $\psi^j \in \dot H^{1/2}_x(\mathbb
R^5)$ such that
\[\lim_{n \to \infty} \big\|v^j(t_n^j + t_n(\lambda_n^j)^2) -
 e^{i(t_n(\lambda_n^j)^{-2} + t_n^j)\Delta}\psi^j\big\|_{\dot H^{1/2}_x(\mathbb R^5)} =0.\]
Thus, together with $(\ref{e317})$ and Lemma $\ref{l210}$ yields
\[
\aligned & \lim_{n \to \infty}\big\langle\,
|\nabla|^{\frac{1}{2}}e^{it_n\Delta}\omega_n^J\, , \,
|\nabla|^{\frac{1}{2}}v^j_n(t_n^j) \,\big\rangle \\
= & \lim_{n \to \infty}\big\langle\,
|\nabla|^{\frac{1}{2}}e^{it_n(\lambda_n^j)^{-2}\Delta}[(g_n^j)^{-1}\omega_n^J]\,
, \,  e^{i(t_n(\lambda_n^j)^{-2} +
t_n^j)\Delta}|\nabla|^{\frac{1}{2}}\psi^j \,\big\rangle \\
= & \lim_{n \to \infty}\big\langle\,
|\nabla|^{\frac{1}{2}}e^{it_n^j\Delta}[(g_n^j)^{-1}\omega_n^J], \,
|\nabla|^{\frac{1}{2}}\psi^j \,\big\rangle \\
 = & \, \,  0.
\endaligned\]
If $\tau^j = -\infty$, then we must have $t_n^j \to -\infty$ as $n
\to \infty$. Indeed, since $t_n(\lambda_n^j)^{-2} \geq 0$ and $
\inf I^j < \infty$, $t_n^j$ can not converges to $+\infty$; if
$t_n^j \equiv 0$, then since $\inf I^j < 0$, we have
$t_n(\lambda_n^j)^{-2} \leq 0$, which contradicts $t_n \in K_n^m
\subset [0, T_{n,j}^+)$. Hence, $\inf I^j = -\infty$. By the
definition of nonlinear profile, $v^j$ scatters backward in time
to $e^{it\Delta}\phi^j$.
\[\lim_{n \to \infty}\big\|v^j(t_n^j + t_n(\lambda_n^j)^2) -
e^{i(t_n(\lambda_n^j)^{-2} + t_n^j)\Delta}\phi^j\big\|_{\dot
H^{1/2}_x(\mathbb R^5)} =0.\] Combining this with $(\ref{e317})$
gives
\[
\aligned & \lim_{n \to \infty}\big\langle\,
|\nabla|^{\frac{1}{2}}e^{it_n\Delta}\omega_n^J \, , \,
|\nabla|^{\frac{1}{2}}v^j_n(t_n^j)  \,\big\rangle \\
= & \lim_{n \to \infty}\big\langle\,
|\nabla|^{\frac{1}{2}}e^{it_n(\lambda_n^j)^{-2}\Delta}[(g_n^j)^{-1}\omega_n^J]\,
, \, e^{i(t_n(\lambda_n^j)^{-2} +
t_n^j)\Delta}|\nabla|^{\frac{1}{2}}\phi^j \,\big\rangle \\
= & \lim_{n \to \infty}\big\langle\,
|\nabla|^{\frac{1}{2}}e^{it_n^j\Delta}[(g_n^j)^{-1}\omega_n^J]\, ,
\,
|\nabla|^{\frac{1}{2}}\phi^j \,\big\rangle \\
 = & \, \, 0.
\endaligned\]
This completes the proof of Step 3.

\vspace{8pt} From $(\ref{e31})$, $(\ref{e313})$, $(\ref{e314})$
\[A_c^2 \geq \limsup_{n \to \infty} \sup_{t \in K_n^m} \|u_n(t)\|^2_{\dot H^{1/2}_x}
\geq \lim_{n \to \infty}\sup_{t\in K_n^m} \Big(\sum_{j=1}^J
\|v_n^j\|^2_{\dot H^{1/2}} + \|\omega_n^J\|^2_{\dot H^{1/2}}\Big).\]

Invoking $(\ref{e312})$ that
\[\limsup_{m \to \infty} \limsup_{n \to \infty} \sup_{t \in K_n^m} \|v_n^{j_0}(t)\|_{\dot H^{1/2}_x} \geq A_c,\]
we conclude that $v_n^j \equiv 0$ for all $j \neq j_0$, and
$\omega_n^{j_0} \to 0 $ in $\dot H^{1/2}_x(\mathbb R^5)$. Thus,
\begin{equation}\label{e318}
u_n(0) = g_ne^{i\tau_n\Delta}\phi + \omega_n
\end{equation}
for some $g_n \in G_{rad}$, $\tau_n \in \mathbb R$, $\phi$,
$\omega_n \in \dot H^{1/2}_x(\mathbb R^5)$ with $\omega_n \to 0$
in $\dot H^{1/2}$. We also have $\tau_n \equiv 0$ or $\tau_n \to
\pm\infty$.

If $\tau_n \equiv 0$, then $u_n(0) \to \phi $ in $\dot H^{1/2}_x$
modulo scaling. This proves Proposition $\ref{p41}$.

If $\tau_n \to \pm\infty$, by time-reversal symmetry, we only
consider $\tau_n \to +\infty$. In this case, by the Strichartz
estimate, we have $\|e^{it\Delta}\phi\|_{S(\mathbb R^+)} <
\infty$. By a change of variables, \[\lim_{n \to
\infty}\|e^{it\Delta}e^{-i\tau_n\Delta}\phi\|_{S(\mathbb R^+)}
=0.\] Taking the group action yields \[ \lim_{n \to
\infty}\|e^{it\Delta}g_ne^{-i\tau_n\Delta}\phi\|_{S(\mathbb R^+)}
=0.\] From $(\ref{e318})$, $(\ref{e22})$, we deduce
\[\lim_{n \to \infty} \|e^{it\Delta}u_n(0)\|_{S(\mathbb R^+)} = 0.\]
Invoking perturbation theorem, we obtain
\[\lim_{n \to \infty}\|u_n\|_{S(\mathbb R^+)} = 0,\]
which contradicts $(\ref{e32})$. This completes the proof of
Proposition $\ref{p41}$.  $\square$ \vspace{8pt}

\subsection{Proof of Theorem 1.3}
{\it Proof.} Suppose Theorem 1.2 failed. Then $A_c <
\frac{\sqrt{6}}{3}\|Q\|_{\dot H^{1/2}_x}$, and  by the definition
of $A_c$, we can find a sequence of solutions $u_n : I_n \times
\mathbb R^5 \mapsto \mathbb C$ to $(1.1)$ with $I_n$ compact,
\begin{equation}\label{e319}
\sup_{n \geq 1} \sup_{t \in I_n}
\big\||\nabla|^{\frac{1}{2}}u_n(t)\big\|_2 = A_c, \quad \lim_{n
\to \infty}\|u_n\|_{S(I_n)} = \infty.
\end{equation}
Then exists $t_n \in I_n$ such that
\begin{equation}\label{e320}
\lim_{n \to \infty}\|u_n\|_{S(-\infty,\; t_n)} = \lim_{n \to
\infty}\|u_n\|_{S(t_n,\;\infty)} = \infty.
\end{equation}
By time-translation symmetry, we set all $t_n =0$. Applying
Proposition $\ref{p41}$, there exists (up to a subsequence) $g_n
\in G_{rad}$ and a function $u_0 \in \dot H^{1/2}_x(\mathbb R^5)$
such that $g_nu_n(0) \to u_0$ in $\dot H^{1/2}_x$. By taking group
action $T_{g_n}$ to the solution $u_n$, we may make $g_n$ be the
identity. Thus $u_n(0) \to u_0$ in $\dot H^{1/2}_x$.

Let $u : I \times \mathbb R^5 \mapsto \mathbb C$ be the
maximal-life span solution to $(1.1)$ with initial data $u(0) =
u_0$. Then, Theorem 1.1 implies $I \subseteq \liminf I_n$ and
\[\lim_{n \to \infty} \sup_{t \in K}\|u_n(t) - u(t)\|_{\dot H^{1/2}_x} = 0\]
for all compact $K \subset I$.

Thus, from $(\ref{e319})$
\begin{equation}\label{e321}
\sup_{t \in I} \|u(t)\|_{\dot H^{1/2}_x} \leq A_c.
\end{equation}

On the other hand, we claim that $u$ blows up both froward and
backward in time. If not, $\|u\|_{S(0,\;\infty)} < \infty$,
$\|u\|_{S(-\infty,\;0)}< \infty$. From perturbation theorem,
$\|u_n\|_{S(0,\;\infty)} < \infty$, $\|u_n\|_{S(-\infty,\;0)} <
\infty$ for $n$ large enough, which contradicts $(\ref{e320})$.

So, by the definition of $A_c$
\[\sup_{t\in I}\|u(t)\|_{\dot H^{1/2}_x} \geq A_c\]
which together with $(\ref{e321})$ yields
\[\sup_{t\in I}\|u(t)\|_{\dot H^{1/2}_x} = A_c.\]

Next, we prove that $u$ is almost periodic modulo scaling. For
arbitrary sequence $\tau_n \in I$, we have
\[\|u\|_{S(-\infty, \;\tau_n)} = \|u\|_{S(\tau_n,\; \infty)} = \infty,\]
since $u$ blows up both forward and backward. From Proposition
$\ref{p41}$, $u(\tau_n)$ has a subsequence which converges in
$\dot H^{1/2}_x(\mathbb R^5)$ modulo scaling. Thus $\{u(t): t\in
I\}$ is precompact in $\dot H^{1/2}_x(\mathbb R^5)$ modulo
$G_{rad}$(Remark 1.2). This completes the proof of the first part
of Theorem 1.3.

 An almost periodic blowup solution which obeys the three scenarios in Theorem 1.3
 can be extracted from the above solution by renormalization  and
 subsequential limits. As we've pointed out, the process is
similar to that in \cite{c13},
 \cite{c14},
 and we refer the readers to these papers for a detailed discussion.

\section{No finite-time blow up}
\setcounter{equation}{0}

 In this section,  we prove
\begin{theorem}
There exists no such maximal life-span solution $u : I \times
\mathbb R^5 \mapsto \mathbb C$ to $(1.1)$  that is  almost periodic
modulo scaling and
\begin{equation}\label{e41}
\sup_{t \in I} \|u(t)\|_{\dot H^{1/2}_x} <
\frac{\sqrt{6}}{3}\|Q\|_{\dot H^{1/2}_x}, \quad \|u\|_{S(I)} =
\infty
\end{equation}
and  either $|\inf I| < \infty$ or $\sup I <\infty$.
\end{theorem}
{\it Proof}. Assume for a contradiction that there existed such a
solution. Without loss of generality, we may assume $\sup I <
\infty$. We claim that
\begin{equation}\label{e42}
\liminf_{t \nearrow \sup I} N(t) = \infty.
\end{equation} If not, we may find a time sequence $t_n \in I$ such
that $t_n \nearrow \sup I$, $\liminf\limits_n N(t_n) < \infty$. For
each $n \geq 1$, define $v_n : I_n \times \mathbb R^5 \mapsto
\mathbb C$ by
\[v_n(t,x) := u(t_n + tN(t_n)^{-2},\; xN(t_n)^{-1}) \]
with $I_n := \{\,t \in \mathbb R:  t_n + tN(t_n)^{-2} \in I\,\}$.
Then $v_n$ is also a solution to $(1.1)$, $\{v_n(0)\}$ is
precompact in $\dot H^{1/2}_x(\mathbb R^5)$. After passing to a
subsequence, we may assume $v_n(0) \to v_0$ in $\dot
H^{1/2}_x(\mathbb R^5)$. Since $\|v_n(0)\|_{\dot H^{1/2}_x} =
\|u(t_n)\|_{\dot H^{1/2}_x}$, $v_0$ is not identically zero.

Let $v$ be the maximal life-span solution to $(1.1)$ with initial
data $v_0$, and maximal life-span $(-T_-,\; T_+)$, $-\infty \leq
T_- < 0 < T_+ \leq \infty$. For  any compact $J \subset (-T_-, \;
T_+)$, from local wellposedness theory, for sufficiently large
$n$, $v_n$ is wellposed on $J$ and $\|v_n\|_{S(J)} < \infty$.
Thus, $u$ is wellposed on the interval $J_n = \{\,t_n +
tN(t_n)^{-2}: t \in J\,\}$ and $\|u\|_{S(J_n)} < \infty$. But
$\liminf_{t \nearrow \sup I} N(t) < \infty$ implies that
$\|u\|_{S}$ is finite beyond $\sup I$, which contradicts the
assumption that $u$ blows up on $I$.

Next, we will prove that for all $R > 0$
\begin{equation}\label{e43}
\limsup_{t \nearrow \sup I} \int_{|x| \leq R} |u(t,x)|^2
\,\mathrm{d} x = 0.
\end{equation}
Let $\eta > 0$, $t \in I$. Using H\"{o}lder's inequality,
Sobolev's embedding theorem, $(\ref{e42})$
\[\aligned
\int_{|x| \leq R}|u(t,x)|^2 \,\mathrm{d} x & \leq \int_{|x|\leq \eta R}
|u(t,x)|^2 \,\mathrm{d} x + \int_{\eta R \leq |x| \leq R}|u(t,x)|^2 \,\mathrm{d} x \\
& \leq \eta R\left(\int |u(t,x)|^{5/2} \,\mathrm{d} x\right)^{4/5} + R
\left(\int_{|x|\geq \eta R} |u(t,x)|^{5/2} \,\mathrm{d} x\right)^{4/5} \\
& \lesssim \eta R\|u(t)\|_{\dot H^{1/2}_x}^2 + R
\left(\int_{|x|\geq \eta R }|u(t,x)|^{5/2} \,\mathrm{d}
x\right)^{4/5}.
\endaligned\]
The first term will go to zero as $\eta$ tends to zero. On the
other hand, from $(\ref{e42})$, almost periodic modulo scaling,
and (\ref{e200}), we have
\[\limsup_{t \nearrow \sup I} \int_{|x|\geq R }|u(t,x)|^{5/2} \,\mathrm{d} x \leq
\limsup_{t \nearrow \sup I} \int_{|x|\geq C(\eta)/N(t)
}|u(t,x)|^{5/2} \,\mathrm{d} x =0. \] Thus $(\ref{e43})$ holds.

We will prove from $(\ref{e43})$ that $u$ is identically zero.

 For $t \in I$,
define
\[M_R(t) := \int_{\mathbb R^5} \phi\big(\frac{|x|}{R}\big)|u(t,x)|^2 \,\mathrm{d} x \]
where $\phi $ is a smooth, radial function with
\[ \phi(r) = \begin{cases}
1, & r \leq 1 \\
0, & r \geq 2.
\end{cases}\]
By $(\ref{e43})$,
\begin{equation}\label{e44}
\limsup_{t \nearrow \sup I}M_R(t) = 0 \quad \textrm{for all } \,\,
R> 0.
\end{equation}
A direct computation involving Plancherel, Hardy's inequality and
$(\ref{e41})$ yields \begin{eqnarray*} \big|\partial_t
M_R(t)\big|&   \lesssim & \int_{\mathbb R^5}
\bigg(\Big|{\frac{x}{R^2}\phi'\big(\frac{|x|}{R}\big)
\bar u} \Big|\bigg)^{\widehat{}}(\xi) |\xi||\hat u| \,\mathrm{d} \xi \\
& \lesssim &
\Big\||\nabla|^{\frac{1}{2}}\big(\frac{x}{R^2}\phi'\big(\frac{|x|}{R}\big)\bar
u
 \big)\Big\|_2 \big\||\xi|^{\frac{1}{2}}\hat u\big\|_2 \\
&  \lesssim_u &
\Big\||\nabla|^{\frac{1}{2}}\big(\frac{x}{R^2}\phi'\big(\frac{|x|}{R}\big)\big)\bar
u\Big\|_2
 +
 \Big\|\frac{x}{R^2}\phi'\big(\frac{|x|}{R}\big)|\nabla|^{\frac{1}{2}}\bar u\Big\|_2
 \\
&  \lesssim_u &
\Big\||x|^{\frac{1}{2}}|\nabla|^{\frac{1}{2}}\big(\frac{x}{R^2}\phi'\big(\frac{|x|}{R}\big)
 \big)\Big\|_{L^\infty}\Big\|\frac{\bar u}{|x|^{1/2}}\Big\|_2 +
 \Big\|\frac{x}{R^2}\phi'\big(\frac{|x|}{R}\big)\Big\|_{L^\infty}\||\nabla|^{\frac{1}{2}}u\|_2.\\
& \lesssim_u &
\Big\||x|^{\frac{1}{2}}|\nabla|^{\frac{1}{2}}\big(\frac{x}{R^2}\phi'\big(\frac{|x|}{R}\big)
 \big)\Big\|_{L^\infty}+ \frac{1}{R} .
\end{eqnarray*}
Furthermore, if $|x| \leq 4R$, then by our chosen of $\phi$
\[\Big\||x|^{\frac{1}{2}}|\nabla|^{\frac{1}{2}}\Big(\frac{x}{R^2}\phi'\big(\frac{|x|}{R}\big)\Big)\Big\|_{L^\infty}
\lesssim \frac{1}{R}.\] If $|x| > 4R$, then using the intrinsic
description of derivatives, we have the following
\[\aligned
\frac{|x|^{\frac{1}{2}}}{R^2}
|\nabla|^{\frac{1}{2}}\Big(x\phi'\big(\frac{|x|}{R}\big) \Big)= &
\, \frac{1}{R^2} \int_{\mathbb R^5}
\frac{|x|^{\frac{1}{2}}\big[x\phi'\big(\frac{|x|}{R}\big)- y
\phi'\big(\frac{|y|}{R}\big)\big]}{|x-y|^{5 + \frac{1}{2}}} \,
\mathrm{d} y \\
= & \frac{1}{R^2} \int_{|x-y| \geq \frac{1}{2}|x|}
\frac{|x|^{\frac{1}{2}}\big[x\phi'\big(\frac{|x|}{R}\big)- y
\phi'\big(\frac{|y|}{R}\big)\big]}{|x-y|^{5 + \frac{1}{2}}} \,
\mathrm{d} y \\
& + \frac{1}{R^2} \int_{|x - y| < \frac{1}{2}|x|}
\frac{|x|^{\frac{1}{2}}\big[x\phi'\big(\frac{|x|}{R}\big)- y
\phi'\big(\frac{|y|}{R}\big)\big]}{|x-y|^{5 + \frac{1}{2}}} \,
\mathrm{d} y.
\endaligned\]
It is easily to check that the first integration has a bound
$R^{-1}$, since $|x- y|\geq \frac{1}{2}|x| \geq 2R$. For the
second one, we have $|y| > |x| - |x - y| > \frac{1}{2}|x| > 2R$,
and by the property of $\phi$, it follows that the integration is
equal to zero.

From the above, we obtain
\[\big|\partial_t M_R(t)\big| \lesssim_u  \frac{1}{R}. \]

 Thus, by the Fundamental Theorem of
Calculus
\[M_R(t_1) \lesssim M_R(t_2) +
\int_{t_2}^{t_1} \partial_t M_R(t) \,\mathrm{d} t \lesssim
M_R(t_2) +  \frac{1}{R}|t_1-t_2|\] for all $t_1, t_2 \in I$ and
$R>0$.

Let $t_2 \nearrow \sup I$  and from $(\ref{e44})$, we obtain
\[M_R(t_1) \lesssim_u  \frac{1}{R} |\sup I - t_1|.\]
Let $R \to \infty$, then we deduce that $M(u(t)) = 0 $ for all $t
\in I$. This implies that $u \equiv 0$, which contradicts
$\|u\|_{S(I)} = \infty$. This completes the proof of Theorem 5.1.

\section{Negative regularity}
\setcounter{equation}{0}

In this section, we prove the following
\begin{theorem}
[Negative regularity in the global case]\label{t61} Let $u$ be a
global radially symmetric solution to $(1.1)$ which is almost
periodic modulo scaling. Suppose also that
\begin{equation}\label{e51}
\sup_{t \in \mathbb R} \|u(t)\|_{\dot H^{1/2}_x} <
\frac{\sqrt{6}}{3}\|Q\|_{\dot H^{1/2}_x}
\end{equation}
and
\begin{equation}\label{e52}
\inf_{t \in \mathbb R} N(t) \gtrsim 1.
\end{equation}
Then, $u \in L_t^\infty \dot H^{-\varepsilon}(\mathbb R \times
\mathbb R^5)$ for some $\varepsilon > 0$. In particular, $u \in
L_t^\infty L_x^2$.
\end{theorem}

In order to prove Theorem $\ref{t61}$, we first establish a
recurrence formula.

Given $\eta > 0$, from Remark 1.1, there exists $N_0=N_0(\eta)$
such that
\begin{equation}\label{e53}
\|u_{\leq N_0}(t)\|_{\dot H^{1/2}_x} \leq \eta.
\end{equation}
Now, define
\[A(N) : = N^{-\frac{3}{4}}\sup_{t \in \mathbb R} \|u_N(t)\|_{L_x^4}\]
for all $N \leq 8N_0$.

Note that by Bernstein's inequality, Sobolev's embedding theorem
\[A(N) \lesssim N^{-\frac{3}{4}}N^{\frac{3}{4}}\|u_N\|_{L_t^\infty L_x^{5/2}} \leq \|u\|_{L_t^\infty \dot H^{1/2}_x} < \infty.\]

Moreover,  $A(N)$ satisfies the following recurrence formula
\begin{lemma}
 For $N \leq 8N_0$
\begin{equation}\label{e54}
A(N) \lesssim_u \left(\frac{N}{N_0}\right)^{\frac{1}{2}} + \eta^2
\sum_{8N \leq N_1 \leq N_0}
\left(\frac{N}{N_1}\right)^{\frac{1}{8}}A(N_1) + \eta^2\sum_{N_1
\leq 8N}\left(\frac{N_1}{N}\right)^{\frac{3}{4}}A(N_1).
\end{equation}
\end{lemma}
{\it Proof}. We only need to prove that for all $t \in \mathbb R$
\[N^{-\frac{3}{4}}\|u_N(t)\|_{L_x^{4}} \lesssim \,\, \textrm{RHS of} \,(\ref{e54}).\]
By the time-translation symmetry, it reduces to prove
\[N^{-\frac{3}{4}}\|u_N(0)\|_{L_x^{4}} \lesssim \,\, \textrm{RHS of} \, (\ref{e54}).\]
By the Duhamel formula $(\ref{e16})$, the triangle, Bernstein's
and the  dispersive inequality, we have
\[\aligned
N^{-\frac{3}{4}}\|u_N(0)\|_{L_x^{4}} & \leq
N^{-\frac{3}{4}}\Big\|\int_0^{N^{-2}} e^{-it\Delta}P_N F(u(t)) \,\mathrm{d}
t
\Big\|_{L_x^4} \\
& \quad + N^{-\frac{3}{4}}\Big\|\int_{N^{-2}}^\infty
e^{-it\Delta}P_N F(u(t)) \,\mathrm{d} t \Big\|_{L_x^4}\\
& \lesssim N^{\frac{1}{2}}\Big\|\int_0^{N^{-2}} e^{-it\Delta}P_N
F(u(t)) \,\mathrm{d} t\Big\|_{L_x^2} \\
& \quad + N^{-\frac{3}{4}}\|P_N F(u)\|_{L_t^\infty L_x^{4/3}}
\int_{N^{-2}}^\infty t^{-\frac{5}{4}} \,\mathrm{d} t \\
& \lesssim N^{-\frac{3}{2}}\|P_N F(u)\|_{L_t^\infty L_x^2} +
N^{-\frac{1}{4}}\|P_N F(u)\|_{L_t^\infty L_x^{4/3}}\\
& \lesssim N^{-\frac{1}{4}}\|P_N F(u)\|_{L_t^\infty L_x^{4/3}}.
\endaligned\]

Decompose $u$ as
\[u := u_{\geq N_0} + u_{\frac{N}{8} \leq \cdot < N_0 } + u_{< \frac{N}{8}},\]
and then make a corresponding expansion of $F(u)$, we obtain terms
constitute $F(u)$ of the following types

1. At least one high frequency, i.e. $|\nabla|^{-2}(u u_{\geq
N_0})u$, or $|\nabla|^{-2}(u^2)u_{\geq N_0}$;

2. Non-high frequency component and at least one lower frequency:
\[|\nabla|^{-2}(u_{< \frac{N}{8}} u_{\leq N_0})u_{\leq N_0}, \quad
|\nabla|^{-2}(u_{\leq N_0}^2)u_{< \frac{N}{8}}; \]

3. All medium components: $|\nabla|^{-2}(u_{\frac{N}{8}\leq \cdot
< N_0}^2)u_{\frac{N}{8} \leq \cdot < N_0}$.

Case 1({\bf{At least one high frequency}}).  Using Bernstein's
inequality, discarding  the projector $P_N$, and then using the
Hardy-Littlewood-Sobolev, H\"{o}lder's and Bernstein's inequality,
Sobolev embedding, we have
\[\aligned  N^{-\frac{1}{4}} \big\|P_N(|\nabla|^{-2}(u
u_{\geq N_0})u)\big\|_{L_t^\infty L_x^{4/3}}  & \lesssim
N^{\frac{1}{2}}\big\||\nabla|^{-2}(u u_{\geq
N_0})u\big\|_{L_t^\infty L_x^{10/9}} \\
& \lesssim_u   N^{\frac{1}{2}}\big\||\nabla|^{-2}(u u_{\geq
N_0})\|_{L_t^\infty
L_x^2}\|u\|_{L_t^\infty L_x^{5/2}} \\
& \lesssim_u N^{\frac{1}{2}}\|u u_{\geq N_0}\|_{L_t^\infty
L_x^{10/9}} \\
& \lesssim_u N^{\frac{1}{2}}\|u\|_{L_t^\infty L_x^{5/2}}\|u_{\geq
N_0}\|_{L_t^\infty L_x^2}\\
 & \lesssim_u N^{\frac{1}{2}}
N_0^{-\frac{1}{2}},
\endaligned\]

\[
\aligned    N^{-\frac{1}{4}} \big\|P_N(|\nabla|^{-2}(u^2) u_{\geq
N_0})\big\|_{L_t^\infty L_x^{4/3}}  &  \lesssim
N^{\frac{1}{2}}\big\||\nabla|^{-2}(u^2)u_{\geq
N_0}\big\|_{L_t^\infty L_x^{10/9}}\\
& \lesssim
N^{\frac{1}{2}}\big\||\nabla|^{-2}(u^2)\big\|_{L_t^\infty
L_x^{5/2}}\|u_{\geq N_0}\|_{L_t^\infty L_x^2}
\\& \lesssim N^{\frac{1}{2}}\|u\|^2_{L_t^\infty L_x^{5/2}}
\|u_{\geq
N_0}\|_{L_t^\infty L_x^2}  \\
& \lesssim_u N^{\frac{1}{2}}N_0^{-\frac{1}{2}};
\endaligned\]

Case 2({\bf Lower frequency components}). By the triangle,
Bernstein's inequality, Sobolev's embedding theorem, H\"{o}lder's
and the Hardy-Littlewood-Sobolev inequality
\[\aligned
 & \qquad N^{-\frac{1}{4}}\big\|P_N(|\nabla|^{-2}(u_{<\frac{N}{8}}u_{\leq
N_0})u_{\leq N_0})\big\|_{L_t^\infty L_x^{4/3}} \\
& \lesssim
N^{-\frac{1}{4}}\big\|P_{>\frac{N}{8}}\big(|\nabla|^{-2}(u_{<
\frac{N}{8}} u_{\leq N_0})\big)u_{\leq N_0}\big\|_{L_t^\infty
L^{4/3}} \\
& \quad + N^{-\frac{1}{4}}\big\||\nabla|^{-2}(u_{<\frac{N}{8}}u_{\leq N_0})P_{>\frac{N}{8}}u_{\leq N_0}\big\|_{L_t^\infty L^{4/3}}\\
& \lesssim N^{-\frac{1}{4}}\big\|P_{>\frac{N}{8}}|\nabla|^{-2}(u_{<\frac{N}{8}}u_{\leq N_0})\big\|_{L_t^\infty L_x^{20/7}}
\|u_{\leq N_0}\|_{L_t^\infty L_x^{5/2}}  \\
& \quad +
N^{-\frac{1}{4}}\big\||\nabla|^{-2}(u_{<\frac{N}{8}}u_{\leq
N_0})\big\|_{L_t^\infty L_x^4}\|P_{>\frac{N}{8}}u_{\leq N_0}\|_{L_t^\infty L_x^2}\\
& \lesssim \eta N^{-\frac{3}{4}}\big\|u_{<\frac{N}{8}}u_{\leq
N_0}\big\|_{L_t^\infty L_x^{20/13}} + N^{-\frac{3}{4}}\big\|u_{<
\frac{N}{8}}u_{\leq N_0}\big\|_{L_t^\infty
L_x^{20/13}}\big\||\nabla|^{\frac{1}{2}}u_{\leq
N_0}\big\|_{L_t^\infty L_x^2} \\ & \lesssim_u \eta^2 \sum_{N_1
\leq \frac{N}{8}} \left(\frac{N_1}{N}\right)^{\frac{3}{4}}A(N_1),
\endaligned\]

\[\aligned
& N^{-\frac{1}{4}}\big\|P_N\big(|\nabla|^{-2}(u_{\leq N_0}^2)u_{<
\frac{N}{8}}\big)\big\|_{L_t^\infty L_x^{4/3}}  \leq
N^{-\frac{1}{4}}\big\|P_{>\frac{N}{4}}|\nabla|^{-2}(u_{\leq
N_0}^2)u_{<
\frac{N}{8}}\big\|_{L_t^\infty L_x^{4/3}} \\
&  \lesssim N^{-\frac{3}{4}}\big\||\nabla|^{-\frac{3}{2}}(u_{\leq
N_0}^2)\big\|_{L_t^\infty L_x^2}\|u_{< \frac{N}{8}}\|_{L_t^\infty
L_x^4}  \lesssim N^{-\frac{3}{4}}\|u_{\leq N_0}^2\|_{L_t^\infty
L_x^{5/4}}\|u_{< \frac{N}{8}}\|_{L_t^\infty L_x^4}\\
&  \lesssim \eta^2 \sum_{N_1 \leq
\frac{N}{8}}\left(\frac{N_1}{N}\right)^{\frac{3}{4}}A(N_1);
\endaligned\]

Case 3({\bf Medium components}). By Bernstein's, the
Hardy-Littlewood-Sobolev, the triangle  and H\"{o}lder's
inequality
\[\aligned
& \quad N^{-\frac{1}{4}}\big\|P_N(|\nabla|^{-2}(u_{\frac{N}{8}
\leq \cdot <
N_0}^2)u_{\frac{N}{8} \leq \cdot <N_0})\big\|_{L_t^\infty L_x^{4/3}}\\
& \lesssim  N^{\frac{1}{8}}\big\||\nabla|^{-2}(u_{\frac{N}{8} \leq
\cdot < N_0}^2) u_{\frac{N}{8} \leq \cdot < N_0}\big\|_{L_t^\infty
L_x^{40/33}}\\
 & \lesssim \sum_{\frac{N}{8} \leq N_1 \leq N_2, \, N_3 \leq N_0}
N^{\frac{1}{8}}\big\||\nabla|^{-2}(u_{N_1}u_{N_2})u_{N_3}\big\|_{L_t^\infty
L_x^{40/33}} \\
& \quad + \sum_{\frac{N}{8} \leq N_3 \leq N_1 \leq N_2 \leq N_0}
N^{\frac{1}{8}}\big\||\nabla|^{-2}(u_{N_1}u_{N_2})u_{N_3}\big\|_{L_t^\infty
L_x^{40/33}}\\
\endaligned\]
\[\aligned
& \lesssim \sum_{\frac{N}{8} \leq N_1 \leq N_2,N_3 \leq N_0}
N^{\frac{1}{8}}\big\||\nabla|^{-2}(u_{N_1}u_{N_2})\big\|_{L_t^\infty
L_x^{40/13}}\big\|u_{N_3}\big\|_{L_t^\infty L_x^2} \\
& \quad + \sum_{\frac{N}{8} \leq N_3 \leq N_1 \leq N_2 \leq N_0}
N^{\frac{1}{8}}
\big\||\nabla|^{-2}(u_{N_1}u_{N_2})\big\|_{L_t^\infty
L_x^{40/23}}\big\|u_{N_3}\big\|_{L_t^\infty L_x^4} \\
& \lesssim_u \sum_{\frac{N}{8} \leq N_1 \leq N_2, N_3 \leq N_0}
N^{\frac{1}{8}} \|u_{N_1}u_{N_2}\|_{L_t^\infty
L_x^{40/29}}N_3^{-\frac{1}{2}} \\
 & \quad + \sum_{\frac{N}{8} \leq N_3 \leq N_1 \leq N_2 \leq N_0}
N^{\frac{1}{8}} \|u_{N_1}u_{N_2}\|_{L_t^\infty L_x^{40/39}}
\|u_{N_3}\|_{L_t^\infty L_x^4}\\
& \lesssim_u \sum_{\frac{N}{8}  \leq N_1 \leq N_2, N_3\leq N_0}
N^{\frac{1}{8}} \|u_{N_1}\|_{L_t^\infty
L_x^4}\|u_{N_2}\|_{L_t^\infty L_x^{40/19}}N_3^{-\frac{1}{2}}\\
& \quad + \sum_{\frac{N}{8} \leq N_3 \leq N_1 \leq N_2 \leq N_0}
N^{\frac{1}{8}} \|u_{N_1}\|_{L_t^\infty
L_x^2}\|u_{N_2}\|_{L_t^\infty L_x^{40/19}}\|u_{N_3}\|_{L_t^\infty
L_x^4}\\
 & \lesssim_u \eta^2\sum_{\frac{N}{8}  \leq N_1 \leq N_2,
N_3\leq N_0} N^{\frac{1}{8}}
N_2^{-\frac{3}{8}}N_3^{-\frac{1}{2}}\|u_{N_1}\|_{L_t^\infty
L_x^4}\\
& \quad + \eta^2 \sum_{\frac{N}{8} \leq N_3 \leq N_1 \leq N_2 \leq
N_0} N^{\frac{1}{8}}
N_1^{-\frac{1}{2}}N_2^{-\frac{3}{8}}\|u_{N_3}\|_{L_t^\infty L_x^4}
\\ &  \lesssim_u \eta^2 \sum_{\frac{N}{8}\leq N_1 \leq
N_0}\left(\frac{N}{N_1}\right)^{\frac{1}{8}}A(N_1).\\
\endaligned\]

This concludes the proof of Lemma 6.1.

\begin{proposition}
Let $u$ be as in Theorem $6.1$. Then
\[u \in L_t^\infty L_x^p \quad \textrm{for} \,\, \frac{22}{9} \leq p < \frac{5}{2},\]
 Furthermore, by the Hardy-Littlewood-Sobolev inequality
 \[|\nabla|^{\frac{1}{2}}F(u) \in L_t^\infty L_x^r \quad \textrm{for} \,\, \frac{110}{101} \leq r < \frac{10}{9}.\]
\end{proposition}
{\it Proof}. Let $N = 8 \cdot 2^{-k}N_0$, applying Lemma 2.1 with
$b_k = (8 \cdot 2^{-k})^{\frac{1}{8}}$, $x_k = A(8 \cdot
2^{-k}N_0)$, we obtain
\[\|u_N\|_{L_t^\infty L_x^4} \lesssim_u N^{7/8+} \quad \textrm{for all } \quad N \leq 8N_0.\]
By the interpolation theorem,  Bernstein's inequality, and
$(\ref{e51})$
\[\aligned
\|u_N\|_{L_t^\infty L^p_x} & \lesssim \|u_N\|_{L_t^\infty
L_x^4}^{2-\frac{4}{p}} \|u_N\|^{\frac{4}{p}-1}_{L_t^\infty
L_x^2}\\
& \lesssim_u N^{\frac{7}{8}(2-\frac{4}{p})+}N^{-\frac{1}{2}(\frac{4}{p}-1)}\\
& \lesssim_u N^{\frac{9}{4}-\frac{11}{2p}+}
\endaligned\]
for all $N \leq 8N_0$.

Thus, using Bernstein's inequality together with $(\ref{e51})$, we
have
\[
\|u\|_{L_t^\infty L_x^p} \leq \|u_{\leq N_0}\|_{L_t^\infty L_x^p}
+ \|u_{> N_0}\|_{L_t^\infty L_x^p} \lesssim_u \sum_{N \leq
N_0}N^{\frac{9}{4}-\frac{11}{2p}+} + \sum_{N >
N_0}N^{2-\frac{5}{p}} \lesssim_u 1.
\]

\begin{proposition}[ Some negative regularity]
Let $u$ be as in Theorem $\ref{t61}$. Assume also that
$|\nabla|^sF(u) \in L_t^\infty L_x^r$ for some $\frac{110}{101}
\leq r < \frac{10}{9}$ and some $0 \leq s \leq \frac{1}{2}$. Then
there exists $s_0 = s_0(r) > 0$ such that $u \in L_t^\infty \dot
H_x^{s-s_0+}$.
\end{proposition}
{\it Proof}. It only needs to prove that
\begin{equation}\label{e55}
\big\||\nabla|^su_N\big\|_{L_t^\infty L_x^2} \lesssim N^{s_0}
\quad \textrm{for all} \quad N>0, \, s_0 := \frac{5}{r}-
\frac{9}{2}.
\end{equation}
In fact, by Bernstein's inequality and $(\ref{e51})$
\[\aligned
\big\||\nabla|^{s-s_0+}u\big\|_{L_t^\infty L_x^2} & \leq
\big\||\nabla|^{s-s_0+}u_{\leq 1}\big\|_{L_t^\infty L_x^2} +
\big\||\nabla|^{s-s_0+}u_{>1}\big\|_{L_t^\infty L_x^2} \\
& \lesssim \sum_{N \leq
1}\big\||\nabla|^{s-s_0+}u_N\big\|_{L_t^\infty
L_x^2} +\sum_{N >1}\big\||\nabla|^{s-s_0+}u_N\big\|_{L_t^\infty L_x^2}\\
& \lesssim_u \sum_{N\leq 1}N^{s_0}N^{-s_0+} + \sum_{N >
1}N^{(s-s_0+)-\frac{1}{2}} \lesssim_u 1.
\endaligned\]

To prove $(\ref{e55})$, by time-translation invariant, we only
need to show that
\[\big\||\nabla|^{s}u_N(0)\big\|_{L_x^2} \lesssim_u N^{s_0} \quad \textrm{for all}\quad N > 0, \, s_0 := \frac{5}{r} - \frac{9}{2} > 0.\]
Using Duhamel formula $(\ref{e16})$ both forward and backward, we
have
\[\aligned
\big\||\nabla|^su_N(0)\big\|_{L_x^2} & = \Big\langle
i\int_0^\infty e^{it\Delta}|\nabla|^sP_NF(u(t)) \,\mathrm{d} t,
-i\int_{-\infty}^0
e^{i\tau\Delta}|\nabla|^sP_NF(u(\tau))\,\mathrm{d} \tau \Big\rangle\\
&\leq \int_0^\infty \int_{-\infty}^0 \Big|\big\langle
e^{it\Delta}|\nabla|^sP_NF(u(t)), \;
e^{i\tau\Delta}|\nabla|^sP_NF(u(\tau))\big\rangle\Big|\,\mathrm{d}
t \,\mathrm{d} \tau.
\endaligned\]
By H\"{o}lder's and the dispersive inequality
\[\aligned
& \quad \Big|\big\langle e^{it\Delta}|\nabla|^sP_NF(u(t)),
e^{i\tau\Delta}|\nabla|^sP_NF(u(\tau))\big\rangle\Big|\\
= & \quad \Big|\big\langle |\nabla|^sP_NF(u(t)),
e^{i(\tau-t)\Delta}|\nabla|^sP_NF(u(\tau))\big\rangle\Big| \\
\leq & \quad
\big\||\nabla|^sP_NF(u(t))\big\|_{L_x^r}\big\|e^{i(\tau-t)\Delta}|\nabla|^sP_NF(u(\tau))\big\|_{L_x^{r'}}\\
\lesssim & \quad |\tau
-t|^{5(\frac{1}{2}-\frac{1}{r})}\big\||\nabla|^sP_NF(u)\big\|_{L_x^r}^2.
\endaligned\]
On the other hand, from Bernstein's inequality
\[\aligned
& \quad \Big|\big\langle e^{it\Delta}|\nabla|^sP_NF(u(t)),
e^{i\tau\Delta}|\nabla|^sP_NF(u(\tau))\big\rangle\Big|\\
\leq & \quad \big\||\nabla|^sP_NF(u)\big\|_{L_x^2}^2\\
\lesssim & \quad
N^{10(\frac{1}{r}-\frac{1}{2})}\big\||\nabla|^sP_NF(u)\big\|_{L_x^r}^2.
\endaligned\]
Thus
\[\aligned
& \,\, \int_0^\infty \int_{-\infty}^0 \Big|\big\langle
e^{it\Delta}|\nabla|^sP_NF(u(t)),
e^{i\tau\Delta}|\nabla|^sP_NF(u(\tau))\big\rangle\Big|\,\mathrm{d}
t \,\mathrm{d} \tau
\\
\lesssim & \,\, \big\||\nabla|^sF(u)\big\|_{L_t^\infty L_x^r}^2
\int_0^\infty \int_{-\infty}^0 \min\{|\tau -
t|^{5(\frac{1}{2}-\frac{1}{r})},
N^{10(\frac{1}{r}-\frac{1}{2})}\}\,\mathrm{d} t \,\mathrm{d} \tau\\
\lesssim  & \,\, \big\||\nabla|^sF(u)\big\|_{L_t^\infty
L_x^r}^2N^{2(\frac{5}{r}-\frac{9}{2})} =  \,\,
\big\||\nabla|^sF(u)\big\|_{L_t^\infty L_x^r}^2N^{2s_0},
\endaligned\]
where we use the fact that $\frac{5}{2}- \frac{5}{r} < -2$.

 With these propositions, we are now
ready to complete the proof of Theorem $\ref{t61}$. First, applying
Proposition 6.2 with $s = \frac{1}{2}$, we obtain $u \in L_t^\infty
\dot H^{\frac{1}{2}-s_0 +}_x$ for some $s_0+ > 0$. By fractional
chain rule and $(\ref{e51})$, we have
$|\nabla|^{\frac{1}{2}-s_0+}F(u) \in L_t^\infty L_x^r$ for some
$\frac{110}{101} \leq r < \frac{10}{9}$. Again using Proposition 6.2
with $s = \frac{1}{2}-s_0+$, we have $u \in L_t^\infty \dot
H^{\frac{1}{2}-2s_0+}_x$. By doing this with finite times, we will
obtain $u \in L_t^\infty \dot H^{-\varepsilon}_x$ for some $0 <
\varepsilon < 2s_0+$. This proves Theorem $\ref{t61}$.

\section{Low-to-high cascade}
\setcounter{equation}{0}

In this section we prove
\begin{theorem}[Absence of cascade]
There can not exist a global solution to $(1.1)$ which is almost
periodic modulo scaling, blows up both forward and backward and is
low-to-high cascade in the sense of Theorem $1.3$.
\end{theorem}
{\it Proof}. We argue by contradiction. Assume there exists such
an $u$. Then, by Theorem 6.1, $u \in L_t^\infty L_x^2$ and
\[0 \leq M(u) = M(u(t)) = \int_{\mathbb R^5} |u(t,x)|^2 \, \mathrm{d} x< \infty \quad \textrm{for all } \quad t \in \mathbb R.\]
Fix $t \in \mathbb R$. Let $\eta > 0$ be sufficiently small. From
$(\ref{e15})$(Remark 1.1)
\[\int_{|\xi| \leq c(\eta)N(t)} |\xi||\hat u(t,\xi)|^2 \,\mathrm{d} \xi \leq \eta.\]
Since $u \in L_t^\infty \dot H^{-\varepsilon}_x(\varepsilon > 0)$,
we see that
\[\int_{|\xi|\leq c(\eta)N(t)}|\xi|^{-2\varepsilon}|\hat u(t,\xi)|^2 \,\mathrm{d} \xi \lesssim 1.\]
Thus, by the interpolation theorem, we obtain
\begin{equation}\label{e61}
\int_{|\xi|\leq c(\eta)N(t)}|\hat u(t,\xi)|^2 \,\mathrm{d} \xi
\lesssim_u \eta^{\frac{2\varepsilon}{1+2\varepsilon}}.
\end{equation}
Meanwhile, it follows from the assumption $(\ref{e51})$ that
\[\aligned
\int_{|\xi|\geq c(\eta) N(t)}|\hat u(t,\xi)|^2 \,\mathrm{d} \xi &
\leq
[c(\eta)N(t)]^{-1}\int |\xi||\hat u(t,\xi)|^2 \,\mathrm{d} \xi\\
& \lesssim_u [c(\eta)N(t)]^{-1}.
\endaligned\]
This together with $(\ref{e61})$ and Plancherel's theorem yields
\[M(u) \lesssim [c(\eta)N(t)]^{-1} + \eta^{\frac{2\varepsilon}{1+2\varepsilon}} \quad \textrm{for all} \quad t \in \mathbb R.\]
As $u$ is a low-to-high cascade solution, there exists $t_n \to
\infty$ such that $N(t_n) \to \infty$. Since $\eta$ is arbitrarily
small, we conclude that $M(u) \equiv 0$. Thus, $u \equiv 0$,
contradicting $\|u\|_{S(\mathbb R)} =0$.

\section{Additional regularity for soliton}
\setcounter{equation}{0}

In order to preclude the final enemy, namely the soliton-like
solution, we need to gain additional regularity to make the
virial-type argument available.
\begin{theorem}
Let $u$ be a global radially symmetric solution to $(1.1)$ that is
almost periodic modulo scaling. Suppose also that $N(t) \equiv 1$
for all $t \in \mathbb R$. Then $u \in L_t^\infty \dot H^s_x$ for
all $s \geq \frac{1}{2}$.
\end{theorem}

To prove Theorem 8.1, we first develop some properties of the
soliton-like solution.

\begin{lemma}[Compactness in $L_x^2$]
Let $u$ be a soliton solution to $(1.1)$ in the sense of Theorem
$1.3$. Then for any $\eta > 0$, there exists $C(\eta) > 0$ such
that
\begin{equation}\label{e71}
\sup_{t \in \mathbb R} \int_{|x| \geq C(\eta)}|u(t,x)|^2 \,
\mathrm{d} x \leq \eta.
\end{equation}
\end{lemma}
{\it Proof}. By negative regularity(Theorem 6.1),
\[\|u_{< N}(t)\|_{L_x^2(|x|\geq R)} \leq
\|u_{< N}(t)\|_{L_x^2} \leq
N^{\varepsilon}\big\||\nabla|^{-\varepsilon}u\big\|_{L_t^\infty
L_x^2} \lesssim_u N^\varepsilon.\] This can be made smaller than
$\eta$ by choosing $N=N(\eta)$ sufficiently small.

To estimate the contribution of high frequency, using Schur's test
lemma \[\big\|\chi_{|x| \geq 2R}
(-\Delta)^{-\frac{1}{2}}|\nabla|^{\frac{1}{2}}P_{\geq N}\chi_{|x|
\leq R}\big\|_{L^2 \to L^2} \lesssim N^{-\frac{1}{2}}\langle
RN\rangle^{-m}.\]

 On the other hand, by Bernstein's inequality
\[\big\|\chi_{|x| \geq 2R} (-\Delta)^{-\frac{1}{2}}|\nabla|^{\frac{1}{2}}P_{\geq N}\chi_{|x| \geq R}\big\|_{L^2 \to L^2}
\lesssim N^{-\frac{1}{2}}.\]
Thus,
\[\aligned
& \,\,  \int_{|x|\geq 2R} |u_{\geq N}(t,x)|^2 \,\mathrm{d} x \\
\lesssim & \,\, \int_{|x|\geq 2R}
\big|(-\Delta)^{-\frac{1}{2}}|\nabla|^{\frac{1}{2}}P_{\geq N}
\chi_{\leq R} |\nabla|^{\frac{1}{2}}u_{\geq N}\big|^2 \,\mathrm{d}
x  \\
& \quad + \int_{|x|\geq 2R}
\big|(-\Delta)^{-\frac{1}{2}}|\nabla|^{\frac{1}{2}}P_{\geq N}
\chi_{\geq
R} |\nabla|^{\frac{1}{2}}u_{\geq N }\big|^2 \,\mathrm{d} x\\
 \lesssim_u & \,\, N^{-1}\langle RN \rangle^{-2m} + N^{-1}\int_{|x|\geq 2R}\big||\nabla|^{\frac{1}{2}}u\big|^2 \,\mathrm{d} x.
\endaligned\]
Choosing $R$ sufficiently large, the first term on the right hand
side can be made smaller than $\eta$. By Definition 1.2, the
second term can also be smaller that $\eta$. Thus, it concludes
$(\ref{e71})$.

\begin{lemma}[Spacetime bounds]
Let $u : I \times \mathbb R^5 \mapsto \mathbb C$ be a maximal
life-span solution to $(1.1)$ which is almost periodic modulo
scaling. Let $J$ be any subinterval of $I$. Then for any
$L^2$-admissible pair $(q, r)$
\begin{equation}\label{s1}
\int_J N(t)^2  \,\mathrm{d} t \lesssim \int_J \Big ( \int_{\mathbb
R^5}\big||\nabla|^{\frac{1}{2}}u(t,x)\big|^r \,\mathrm{d}
x\Big)^{q/r} \,\mathrm{d} t \lesssim 1 + \int_J N(t)^2
\,\mathrm{d} t
\end{equation}
\end{lemma}
{ \it Proof}.  As noted, the proof can be found in \cite{c13},
\cite{c14}. For the sake of convenience, we give a self-contained
argument using the ideas in them.

We first prove the second inequality. Let $\eta > 0$ be chosen
later, divide $J$ into subintervals $I_j$ such that on each $I_j$
\[
\int_{I_j} N(t)^2 \, \mathrm{d} t \leq \eta.
\]
By pigeonhole principle, there are at most $m \leq \eta^{-1}
\times \big( 1 + \int_J N(t)^2 \, \mathrm{d} t\big)$ subintervals.
For each $j$,  choose $t_j$ such that
\begin{equation}\label{s2}
N(t_j)^2|I_j| \leq 2\eta.
\end{equation}
By Strichartz's estimate, the Hardy-Littlewood-Sobolev, and
H\"{o}lder's, Sobolev's inequality, we have on $I_j \times \mathbb
R^5$ that
\[
\aligned \big\| |\nabla|^{\frac{1}{2}}u \big\|_{L_t^q L_x^r} \leq
& \quad \big\| e^{i(t-t_j)\Delta}|\nabla|^{\frac{1}{2}}u(t_j)
\big\|_{L_t^q L_x^r} \\
& \qquad + \Big\|\int_{t_j}^t
e^{i(t-\tau)\Delta}|\nabla|^{\frac{1}{2}}F(u(\tau)) \,
\mathrm{d}\tau \Big\|_{L_t^q L_x^r} \\
\lesssim & \quad \big\| |\nabla|^{\frac{1}{2}}u_{\geq
N_0}(t_j)\big\|_2 + \big\|e^{i(t
-t_j)\Delta}|\nabla|^{\frac{1}{2}}u_{\leq N_0}(t_j)
\big\|_{L_t^q L_x^r} \\
& \qquad + \big\| |\nabla|^{\frac{1}{2}}F(u)\big\|_{L_t^{\tilde
q'}L_x^{\tilde r'}} \\
\lesssim & \quad \big\| |\nabla|^{\frac{1}{2}}u_{\geq
N_0}(t_j)\big\|_2 + |I_j|^{1/q}N_0^{2/q}\big\|
|\nabla|^{\frac{1}{2}}u_{< N_0}
\big\|_{L_t^\infty L_x^2}\\
& \qquad + \big\||\nabla|^{\frac{1}{2}}u \big\|^3_{L_t^q L_x^r},
\endaligned\]
where $\tilde q' = q/3$, $\tilde r' = (15-3r)/5r$.  From the
definition of almost periodic modulo scaling, choosing $N_0$ as a
large multiple of $N(t_j)$, then the first term on the right hand
side can be made as small as we wish. Invoking $(\ref{s2})$ and
choosing $\eta$ sufficiently small, the second term can also be
made sufficiently small. Thus, by bootstrap argument, we obtain
\[\int_{I_j} \Big( \int_{\mathbb R^5}\big| |\nabla|^{\frac{1}{2}}u(t,x)\big|^r \, \mathrm{d} x\Big)^{q/r} \, \mathrm{d}t \leq \eta. \]
Recalling the bound on subinterval number, we have
\[\int_{J} \Big( \int_{\mathbb R^5}\big| |\nabla|^{\frac{1}{2}}u(t,x)\big|^r \, \mathrm{d} x\Big)^{q/r} \, \mathrm{d}t \leq
1 + \int_J N(t)^2 \, \mathrm{d} t.\] For the first inequality,
note that by Definition 1.2, we must have
\[\int_{|x| \leq C(\eta)N(t)^{-1}} \big| |\nabla|^{\frac{1}{2}}u(t,x)\big|^2 \, \mathrm{d} x \gtrsim_u 1.\]
Using H\"{o}lder's inequality
\[\Big(\int_{\mathbb R^5} \big| |\nabla|^{\frac{1}{2}}u(t,x)\big|^r \, \mathrm{d} x\Big)^{1/r}
\gtrsim \Big(\int_{|x|\leq
C(\eta)N(t)^{-1}}\big||\nabla|^{\frac{1}{2}}u(t,x)\big|^2 \,
\mathrm{d} x\Big)^{1/2} N(t)^{2/q} \gtrsim_u N(t)^{2/q}.\]
Integrating the above inequality on $J$, we have
\[\int_J\Big(\int_{\mathbb R^5} \big| |\nabla|^{\frac{1}{2}}u(t,x)\big|^r \, \mathrm{d}x\Big)^{q/r}\, \mathrm{d} t
\gtrsim_u \int_J N(t)^2 \, \mathrm{d} t .\]

\begin{remark} We have for all $\dot H^{1/2}$-admissible pairs $(q,\, r)$
that
\[\int_J N(t)^2  \,\mathrm{d} t \lesssim \|u\|_{L_t^q L_x^r(J \times \mathbb R^5)}^q \lesssim 1 + \int_J N(t)^2
\,\mathrm{d} t.\]
Indeed,
\[
\aligned \|u \|_{L_t^q L_x^r} \leq & \quad \big\|
e^{i(t-t_j)\Delta}|\nabla|^{\frac{1}{2}}u(t_j)
\big\|_{L_t^q L_x^r} \\
& \qquad + \Big\|\int_{t_j}^t
e^{i(t-\tau)\Delta}|\nabla|^{\frac{1}{2}}F(u(\tau)) \,
\mathrm{d}\tau \Big\|_{L_t^q L_x^r} \\
\lesssim & \quad \big\| |\nabla|^{\frac{1}{2}}u_{\geq
N_0}(t_j)\big\|_2 + |I_j|^{1/q}N_0^{2/q}\big\|
|\nabla|^{\frac{1}{2}}u_{< N_0}
\big\|_{L_t^\infty L_x^2}\\
& \qquad + \|u\|_{L_t^q L_x^r}^2\big\||\nabla|^{\frac{1}{2}}u
\big\|_{L_t^{q_1} L_x^{r_1}},
\endaligned\]
where $(q_1, \, r_1)$ is an $L^2$-admissible pair. Using the same
argument as that in proving $(\ref{s1})$, we easily get the
bounds.
\end{remark}
 Due to this proposition, we could obtain some
local estimates for the soliton-like solution. Specifically, we
have for $L^2$-admissible pair $(q, r)$ and $\dot
H^{1/2}$-admissible pair $(\tilde q, \tilde r)$ that
\begin{equation}\label{e73}
\|u\|_{L_t^{\tilde q}L_x^{\tilde r}(J \times \mathbb R^5)}
\lesssim_u \langle |J| \rangle^{\frac{1}{\tilde q}}, \quad
\big\||\nabla|^{\frac{1}{2}} u\big\|_{L_t^qL_x^r(J \times \mathbb
R^5)} \lesssim_u \langle |J| \rangle^{\frac{1}{q}}.
\end{equation}
By the Hardy-Littlewood-Sobolev inequality and the interpolation
\begin{eqnarray}\label{e74}
 \|F(u)\|_{L_t^2L_x^{10/7}}  & \leq \|(|\cdot|^{-3} \ast
|u|^2)\|_{L_t^2L_x^{10/3}}\|u\|_{L_t^\infty L_x^{5/2}} \lesssim_u
\|u\|_{L_t^4 L_x^{20/7}}^2 \nonumber\\
& \lesssim_u \|u\|_{L_t^4
L_x^{10/3}}\big\||\nabla|^{\frac{1}{2}}u\big\|_{L_t^4 L_x^{5/2}}
\lesssim_u \langle |J| \rangle^{\frac{1}{2}}.
\end{eqnarray}
 By the weighted Strichartz estimate
\begin{equation}\label{e75}
\big\||x|^2 u\big\|_{L^4_tL_x^\infty} \lesssim_u \langle |J|
\rangle^{\frac{1}{2}}.
\end{equation}
 From Definition 1.2
\begin{equation}\label{e76}
\lim_{N \to \infty}\|u_{\geq N}\|_{L_t^\infty \dot
H^{1/2}_x(\mathbb R \times \mathbb R^5)} = 0 .\end{equation}

Now, define
\begin{equation}\label{e77}
G(N) := \|u_{\geq N}\|_{L_t^\infty \dot H^{1/2}_x(\mathbb R \times
\mathbb R^5)}
\end{equation}
Note that
\begin{equation}\label{e78}
\lim_{N \to \infty} G(N) = 0.
\end{equation}

To prove Theorem 8.1, it suffices to prove that $G(N) \lesssim_u
N^{-s}$ holds for any $s > 0$ and any sufficiently large $N$,
since we consequently have $\|u_N\|_{L_t^\infty \dot H_x^{s +
1/2}} \lesssim N^s\|u_N\|_{L_t^\infty \dot H^{1/2}_x} \lesssim_u 1
$. This will be achieved by iterating the following proposition
with sufficiently small $\eta$.
\begin{proposition}
Let $u$ be as in Theorem $8.1$. Let $\eta > 0$ be sufficiently
small. Then, for sufficiently large $N=N(\eta, u)$, we have
\begin{equation}\label{e79}
G(N) \lesssim_u \eta G\big(\frac{N}{16}\big).
\end{equation}
\end{proposition}
 To prove the proposition, it suffices to prove
\begin{equation}\label{e710}
\|u_{\geq N}(t)\|_{\dot H^{1/2}_x} \lesssim_u \eta
G\left(\frac{N}{16}\right)
\end{equation}
for all $t \in \mathbb R$ and all $N$ sufficiently large. By
time-translation invariant, we may set $t = 0$. Using Duhamel
formula $(\ref{e16})$ and the in/out decomposition
\begin{eqnarray}\label{e711}
|\nabla|^{\frac{1}{2}}u_{\geq N}(0) & = & (P^+ +
P^-)|\nabla|^{\frac{1}{2}}u_{\geq N}(0) \nonumber\\
& = & \lim_{T \to \infty}i \int_0^T P^+e^{-it\Delta}P_{\geq
N}|\nabla|^{\frac{1}{2}}F(u(t))\,\mathrm{d} t \nonumber \\
& & - \lim_{T \to \infty} \int_{-T}^0 P^-e^{-it\Delta}P_{\geq
N}|\nabla|^{\frac{1}{2}}F(u(t))\,\mathrm{d} t
\end{eqnarray}
as weak limits in $L_x^2$. Using the property of weak closedness
for unit ball, namely \[f_T \rightharpoonup f \quad
\Longrightarrow \quad \|f\| \leq \liminf_{T} \|f_T\|,\] we are
reduced to proving that RHS of $(\ref{e711})$
 $\lesssim_u$ RHS of $(\ref{e710})$.

Note that $P^{\pm}$ are singular at $x=0$; to get around this, we
introduce the cutoff $\psi_N(x) := \psi(N|x|)$, where $\psi$ is
the characteristic function of $[1, \infty)$. As the short times
and large times will be treated differently, we rewrite
$(\ref{e711})$ as
\begin{eqnarray*}
|\nabla|^{\frac{1}{2}}u_{\geq N}(0) & = & [\psi_N(x) + (1 -
\psi_N(x))]|\nabla|^{\frac{1}{2}}u_{\geq N}(0) \nonumber \\
& = &\lim_{T \to \infty}\int_0^T \psi_N(x)P^+e^{-it\Delta}P_{\geq
N}|\nabla|^{\frac{1}{2}}F(u(t)) \,\mathrm{d} t \nonumber \\
& & - \lim_{T \to \infty}i \int_{-T}^0
\psi_N(x)P^-e^{-it\Delta}P_{\geq
N}|\nabla|^{\frac{1}{2}}F(u(t)) \,\mathrm{d} t \nonumber \\
& & + \lim_{T \to \infty}i \int_0^T (1 -
\psi_N(x))e^{-it\Delta}P_{\geq
N}|\nabla|^{\frac{1}{2}}F(u(t)) \,\mathrm{d} t \nonumber, \\
\end{eqnarray*}
\begin{eqnarray}\label{e712}
|\nabla|^{\frac{1}{2}}u_{\geq N}(0)& = & i\int_0^\delta
\psi_N(x)P^+e^{-it\Delta}P_{\geq
N}|\nabla|^{\frac{1}{2}}F(u(t)) \,\mathrm{d} t \nonumber \\
& & -i\int_{-\delta}^0 \psi_N(x)P^-e^{-it\Delta}P_{\geq
N}|\nabla|^{\frac{1}{2}}F(u(t)) \,\mathrm{d} t \nonumber \\
& & + i\int_0^\delta (1 - \psi_N(x))e^{-it\Delta}P_{\geq
N}|\nabla|^{\frac{1}{2}}F(u(t)) \,\mathrm{d} t \nonumber \\
& & + \lim_{T \to \infty}\sum_{M \geq N}i\int_{\delta}^T
\int_{\mathbb R^5}\psi_N[P^+_Me^{-it\Delta}](x,y) \tilde P_M|\nabla|^{\frac{1}{2}}F(u(t)) \,\mathrm{d} y \,\mathrm{d} t \nonumber \\
& &  - \lim_{T \to \infty}\sum_{M \geq N}i\int_{-T}^{-\delta}
\int_{\mathbb R^5}\psi_N[P^-_Me^{-it\Delta}](x,y) \tilde P_M|\nabla|^{\frac{1}{2}}F(u(t)) \,\mathrm{d} y \,\mathrm{d} t \nonumber \\
&  & + \lim_{T \to \infty}  \sum_{M \geq N}i \int_{\delta}^T
\int_{\mathbb R^5}(1 - \psi_N)[\tilde P_Me^{-it\Delta}](x,y)  P_M|\nabla|^{\frac{1}{2}}F(u(t)) \,\mathrm{d} y \,\mathrm{d} t \nonumber \\
& := & I_1 - I_2 + I_3 + I_4 - I_5 + I_6.
\end{eqnarray}
Note that we used the identity
\[P_{\geq N} = \sum_{M \geq N} P_M \tilde P_M,\]
where $\tilde P_M := P_{M/2} + P_M + P_{2M}$.

For integrals over short times, namely $I_1$, $I_2$, $I_3$, we
have the following estimate, that is
\begin{lemma}[Local estimate]
For any sufficiently small $\eta > 0$, there exists $\delta =
\delta(u, \eta) > 0$ such that
\[
\Big\| \int_0^\delta e^{-it\Delta}P_{\geq
N}|\nabla|^{\frac{1}{2}}F(u(t)) \,\mathrm{d} t\Big\|_{L_x^2}
\lesssim_u \eta G\left(\frac{N}{8}\right)
\]for sufficiently large $N$ depending on $u$ and $\eta$. An
analogous estimate holds for integration over $[-\delta, 0]$ and
after pre-multiplication by $P^{\pm}$.
\end{lemma}
{\it Proof}. By Strichartz's estimate, it only needs to prove
\[\big\||\nabla|^{\frac{1}{2}}P_{\geq N}F(u)\big\|_{L_t^2L_x^{10/7}(J \times \mathbb R^5)} \lesssim_u \eta G\left(\frac{N}{8}\right)\]
for any time interval $J$ with $|J| \leq \delta$.

From $(\ref{e78})$, for any $\eta > 0$, there exists $N_0 =
N_0(u,\eta)$ such that \begin{equation}\label{e713}\|u_{\geq
N_0}\|_{L_t^\infty \dot H^{1/2}_x} \leq \eta.
\end{equation}
 Let $N
\geq 8N_0$. Decompose $u$ as
\[u := u_{\geq \frac{N}{8}} + u_{N_0 \leq \cdot < \frac{N}{8}} + u_{< N_0},\]
and make a corresponding expansion of $P_{\geq N}F(u)$. Note that
any term in the resulting expansion does not contain $u_{\geq
\frac{N}{8}}$ vanishes.

We first consider a term with two factors of the form $u_{< N_0}$.
Using H\"{o}lder's inequality, the fractional Leibniz rule, the
Hardy-Littlewood-Sobolev, and Bernstein's inequality
\[\aligned
& \quad \big\||\nabla|^{\frac{1}{2}}(|\nabla|^{-2}(u_{<
N_0}^2)u_{\geq
\frac{N}{8}})\big\|_{L_t^2L_x^{10/7}{(J \times \mathbb R^5)}} \\
& \leq \quad  \big\||\nabla|^{-2}(u_{<
N_0}^2)\big\|_{L_t^2L_x^5{(J \times \mathbb
R^5)}}\big\||\nabla|^{\frac{1}{2}}u_{\geq
\frac{N}{8}}\big\|_{L_t^\infty L_x^2} \\
&  \qquad + \big\||\nabla|^{-\frac{3}{2}}(u_{<
N_0}^2)\big\|_{L_t^2L_x^{10/3}{(J \times \mathbb R^5)}}\|u_{\geq
\frac{N}{8}}\|_{L_t^\infty L_x^{5/2}} \\
& \lesssim \quad \|u_{< N_0}^2\|_{L_t^2L_x^{5/3}(J \times \mathbb
R^5)}G\big(\frac{N}{8}\big) + \|u_{< N_0}^2\|_{L_t^2L_x^{5/3}(J
\times \mathbb R^5)}G\big(\frac{N}{8}\big)\\
& \lesssim_u  \quad |J|^{\frac{1}{2}}N_0 G\big(\frac{N}{8}\big),
\endaligned\]
and
\[\aligned
 & \quad \big\||\nabla|^{\frac{1}{2}}(|\nabla|^{-2}(u_{< N_0}u_{\geq
\frac{N}{8}})u_{< N_0})\big\|_{L_t^2L_x^{10/7}{(J \times \mathbb R^5)}} \\
 & \leq  \quad \big\||\nabla|^{-2}(u_{< N_0}u_{\geq
\frac{N}{8}})\big\|_{L_t^4L_x^{10/3}}\big\||\nabla|^{\frac{1}{2}}u_{<
N_0}\big\|_{L_t^4L_x^{5/2}} \\
&  \qquad + \big\||\nabla|^{-\frac{3}{2}}(u_{< N_0}u_{\geq
\frac{N}{8}})\big\|_{L_t^4L_x^{5/2}}\|u_{<
N_0}\|_{L_t^4L_x^{10/3}}\\
& \lesssim_u \quad \|u_{< N_0}u_{\geq
\frac{N}{8}}\|_{L_t^4L_x^{10/7}}|J|^{\frac{1}{4}}N_0^{\frac{1}{2}}
+ \|u_{< N_0}u_{\geq
\frac{N}{8}}\|_{L_t^4L_x^{10/7}}|J|^{\frac{1}{4}}N_0^{\frac{1}{2}}\\
& \lesssim_u \quad \|u_{< N_0}\|_{L_t^4L_x^{10/3}}\|u_{\geq
\frac{N}{8}}\|_{L_t^\infty
L_x^{5/2}}|J|^{\frac{1}{4}}N_0^{\frac{1}{2}}\\
 & \lesssim_u \quad
|J|^{\frac{1}{2}}N_0 G\big(\frac{N}{8}\big) .
\endaligned\]
Choosing $\delta$ sufficiently small depending on $\eta$ and
$N_0$, we see they are acceptable.

Now, we have to estimate those components of $P_{\geq N}F(u)$
which involve $u_{\geq \frac{N}{8}}$ and at least one of the other
terms is not $u_{< N_0}$. Using H\"{o}lder's inequality, the
fractional Leibniz rule, the Hardy-Littlewood-Sobolev, Bernstein's
inequality, $ (\ref{e73})$, $(\ref{e713})$,
\[\aligned
& \quad \big\||\nabla|^{\frac{1}{2}}(|\nabla|^{-2}(u_{\geq
N_0}u_{\geq
\frac{N}{8}})u)\big\|_{L_t^2L_x^{10/7}(J \times \mathbb R^5)}\\
& \lesssim \quad \big\||\nabla|^{-\frac{3}{2}}(u_{\geq N_0}u_{\geq
\frac{N}{8}})\big\|_{L_t^4L_x^{5/2}(J\times \mathbb
R^5)}\|u\|_{L_t^4L_x^{10/3}(J\times \mathbb R^5)} \\
& \qquad + \big\||\nabla|^{-2}(u_{\geq N_0}u_{\geq
\frac{N}{8}})\big\|_{L_t^\infty L_x^{5}(J\times \mathbb
R^5)}\big\||\nabla|^{\frac{1}{2}}u\big\|_{L_t^2 L_x^2}\\
& \lesssim_u \quad \|u_{\geq N_0}u_{\geq
\frac{N}{8}}\|_{L_t^4L_x^{10/7}(J\times \mathbb R^5)}\langle
|J|\rangle^{\frac{1}{4}} + |J|^{\frac{1}{2}}\|u_{\geq
N_0}u_{\geq\frac{N}{8}}\|_{L_t^\infty
L_x^{5/3}(J \times \mathbb R^5)} \\
& \lesssim_u \quad \|u_{\geq \frac{N}{8}}\|_{L_t^\infty
L_x^{5/2}}\|u_{\geq N_0}\|_{L_t^4L_x^{10/3}(J \times \mathbb
R^5)}\langle |J| \rangle^{\frac{1}{4}} +
|J|^{\frac{1}{2}}\|u_{\geq \frac{N}{8}}\|_{L_t^\infty
L_x^{5/2}}\|u_{\geq N_0}\|_{L_t^\infty L_x^{5}(J\times \mathbb
R^5)} \\
& \lesssim_u  \quad \langle |J| \rangle^{\frac{1}{4}}\|u_{\geq
N_0}\|_{L_t^4L_x^{10/3}(J \times \mathbb
R^5)}G\big(\frac{N}{8}\big) +
\eta|J|^{\frac{1}{2}}N_0G\big(\frac{N}{8}\big)
\endaligned\]
By $(\ref{e73}) $
\[\|u_{\geq N_0}\|_{L_t^2 L_x^5(J\times \mathbb R^5)} \lesssim \langle |J| \rangle^{\frac{1}{2}}.\]
Hence, interpolating with $(\ref{e713})$, we have
\[\|u_{\geq N_0}\|_{L_t^4L_x^{10/3}(J \times \mathbb
R^5)} \lesssim \eta^{\frac{1}{2}}\langle |J|
\rangle^{\frac{1}{4}}.\] Thus, we obtain
\[\big\||\nabla|^{\frac{1}{2}}(|\nabla|^{-2}(u_{\geq N_0}u_{\geq
\frac{N}{8}})u)\big\|_{L_t^2L_x^{10/7}(J \times \mathbb R^5)}
\lesssim_u \eta^{\frac{1}{2}}\langle
|J|\rangle^{\frac{1}{2}}G\big(\frac{N}{8}\big) + \eta
|J|^{\frac{1}{2}}N_0G\big(\frac{N}{8}\big),\] which is acceptable.

In the same manner, we estimate
\[\aligned
& \quad \big\||\nabla|^{\frac{1}{2}}(|\nabla|^{-2}(u_{\geq
\frac{N}{8}}u)u_{\geq N_0})\big\|_{L_t^2L_x^{10/7}(J \times
\mathbb
R^5)}\\
& \lesssim   \big\||\nabla|^{-\frac{3}{2}}(u_{\geq
\frac{N}{8}}u)\big\|_{L_t^2 L_x^{10/3}(J \times \mathbb R^5)}
\|u_{\geq N_0}\|_{L_t^\infty L_x^{5/2}} \\
& \qquad + \big\||\nabla|^{-2}(u_{\geq
\frac{N}{8}}u)\big\|_{L_t^\infty L_x^{5}(J \times \mathbb
R^5)}\big\||\nabla|^{\frac{1}{2}}u_{\geq N_0}\big\|_{L_t^2 L_x^2}\\
& \lesssim_u  \eta \|u\|_{L_t^2L_x^5(J\times \mathbb R^5)}
\|u_{\geq \frac{N}{8}}\|_{L_t^\infty L_x^{5/2}}  \lesssim_u
\eta\langle |J|\rangle^{\frac{1}{2}}G\big(\frac{N}{8}\big).
\endaligned\]
Another term
$\big\||\nabla|^{\frac{1}{2}}\big(|\nabla|^{-2}(u_{\geq
N_0}u)u_{\geq \frac{N}{8}}\big)\big\|_{L_t^2L_x^{10/7}(J \times
\mathbb R^5)}$ can be estimated similarly. This concludes the
proof of Lemma 8.3.

We now turn our attention to $I_4, \, I_5, \, I_6$, namely the
integrations over large times: $|t| \geq \delta$. Making use of
the properties of the kernels $P_Me^{-it\Delta}$,
$P_M^{\pm}e^{-it\Delta}$(see Lemma 2.4, Lemma 2.7), we break the
regions of $(t,y)$ integration into two pieces: $|y| \gtrsim M|t|$
and $|y| \ll M|t|$. when $|x|\geq N^{-1} $, we use the kernel
$P_M^{\pm}e^{-it\Delta}$; in this case $|y|- |x| \thicksim M|t|$
implies $|y| \gtrsim M|t|$ for $|t| \geq \delta \geq N^{-2}$. When
$|x| \leq N^{-1} $, we use $P_Me^{-it\Delta}$; in this case $|y-x|
\thicksim M|t|$ implies $|y| \gtrsim M|t|$ for $ |t| \geq \delta
\geq N^{-2}$. The condition $\delta \geq N^{-2}$ can be satisfied
under our statement $N$ sufficiently large depending on $u$ and
$\eta$.

 Define $\chi_k$ as the
characteristic function of the set
\[\{\,(t,y): 2^k\delta \leq |t| \leq 2^{k+1}\delta, |y| \gtrsim M|t|\,\}.\]
Then we have the following estimate
\begin{lemma}[Main contribution]
Let $\eta > 0$ be a small number and $\delta$ be as in Lemma
$8.3$. Then
\begin{equation}\label{e714}
\sum_{M\geq N}\sum_{k=0}^\infty
\Big\|\int_{2^k\delta}^{2^{k+1}\delta}\int_{\mathbb R^5}
[P_Me^{-it\Delta}](x,y)\chi_k(t,y)[\tilde P_M
|\nabla|^{\frac{1}{2}}F(u(t))](y)\,\mathrm{d} y \,\mathrm{d} t
\Big\|_{L_x^2} \lesssim_u \eta G\big(\frac{N}{16}\big)
\end{equation}
for all $N$ sufficiently large depending on $u$ and $\eta$. An
analogous estimate holds for integration over $[-2^{k+1}\delta,
-2^k\delta]$ and with $P_M$ replaced by $P_M^{\pm}$.
\end{lemma}
{\it Proof}. By Strichartz's estimates
\begin{eqnarray*}
& &  \Big\|\int_{2^k\delta}^{2^{k+1}\delta}\int_{\mathbb
R^5}[P_Me^{-it\Delta}](x,y)\chi_k(t,y)[\tilde P_M
|\nabla|^{\frac{1}{2}}F(u(t))](y) \,\mathrm{d} y \,\mathrm{d} t\Big\|_{L_x^2}\\
& & \lesssim \big\|\chi_k \tilde
P_M(|\nabla|^{\frac{1}{2}}F(u))\big\|_{L_t^2L_y^{10/7}([2^k\delta,
\, 2^{k+1}\delta]\times \mathbb R^5)}
\end{eqnarray*}
Using the fractional Leibniz rule, we turn to estimate
\[ {\rm II_1} =
\big\|\chi_k\tilde
P_M(|\nabla|^{-\frac{3}{2}}(|u|^2)u)\big\|_{L_t^2L_x^{10/7}([2^k\delta,\,
2^{k+1}\delta]\times \mathbb R^5)},\] \[ {\rm II_2}=
\big\|\chi_k\tilde
P_M(|\nabla|^{-2}(|u|^2)|\nabla|^{\frac{1}{2}}u)\big\|_{L_t^2L_x^{10/7}([2^k\delta,
\,2^{k+1}\delta]\times \mathbb R^5)}.\]

We only estimate ${\rm II_1}$, since ${\rm II_2}$ can be treated
similarly, using the fact that $u \in L_t^\infty H^{1/2}$.

Write $u$ as $u := u_{\leq {\frac{M}{16}}} + u_{> \frac{M}{16}}$.
In what follows, all spacetime norms are taken on the slab
$[2^k\delta,\; 2^{k+1}\delta] \times \mathbb R^5$, unless noted
otherwise. Using the support property of $\tilde P_M$, ${\rm
II_{11}}$ can be controlled by
\begin{eqnarray*}
 {\rm II_1}  & \lesssim & \big\| \chi_k |\nabla|^{-\frac{3}{2}}(u^2)u_{> \frac{M}{16}}
 \big\|_{L_t^2 L_x^{10/7}} + \big\| \chi_k |\nabla|^{-\frac{3}{2}}(u u_{> \frac{M}{16}})u_{\leq \frac{M}{16}}
 \big\|_{L_t^2 L_x^{10/7}} \\
& := & {\rm II_{11} + II_{12}}.
\end{eqnarray*}

Using H\"{o}lder's inequality, and $(\ref{e75})$, we have
\begin{eqnarray*}
& & \big\| \chi_k |\nabla|^{-\frac{3}{2}}(u^2)u_{> \frac{M}{16}}
 \big\|_{L_t^2 L_x^{10/7}} \leq \big\|u_{> \frac{M}{16}}\big\|_{L_t^\infty
L_x^{5/2}}\big\|\chi_k|\nabla|^{-\frac{3}{2}}(u^2)\big\|_{L_t^2L_x^{10/3}}\\
 & \lesssim & G\big(\frac{M}{16}\big)\left( \Big\|\chi_k \int_{|x-y|
\geq \frac{|y|}{2}}\frac{|u(x)|^2}{|x-y|^{7/2}}\,\mathrm{d}
x\Big\|_{L_t^2L_y^{10/3}} + \Big\|\chi_k \int_{|x-y| <
\frac{|y|}{2}}\frac{|u(x)|^2}{|x-y|^{7/2}}\,\mathrm{d}
x\Big\|_{L_t^2L_y^{10/3}}\right) \\
& \lesssim & G\big(\frac{M}{16}\big) \left(\|\chi_k
|y|^{-\frac{7}{2}}\|_{L_t^2L_y^{10/3}}\|u\|_{L_t^\infty L_x^2} +
\Big\|\chi_k |y|^{-\frac{16}{5}}\int_{|x-y| <
\frac{|y|}{2}}\frac{|y|^{16/5}|u|^2}{|x-y|^{7/2}} \,\mathrm{d} x
\Big\|_{L_t^2L_y^{10/3}} \right)\\
 & \lesssim_u &
G\big(\frac{M}{16}\big)\left(M^{-2}(2^k\delta)^{-\frac{3}{2}} +
 \Big\|\chi_k |y|^{-\frac{16}{5}}\big\|1_{\leq \frac{|y|}{2}}|\cdot|^{-\frac{7}{2}}\big\|_{L_x^{5/4}}
 \big\||y|^2u\|^{\frac{8}{5}}_{L_x^\infty}\|u\|
 _{L_x^2}^{\frac{2}{5}} \Big\|_{L_t^2L_y^{10/3}}\right)\\
& \lesssim_u &
G\big(\frac{M}{16}\big)\left(M^{-2}(2^k\delta)^{-\frac{3}{2}} +
  \big\|\chi_k
|y|^{-\frac{27}{10}}\big\|_{L_t^{10}L_y^{10/3}}\big\||y|^2u\big\|_{L_t^4L_x^\infty}^{\frac{8}{5}}\right)\\
& \lesssim_u &
G\big(\frac{M}{16}\big)\left(M^{-2}(2^k\delta)^{-\frac{3}{2}} +
M^{-\frac{6}{5}}(2^k\delta)^{-\frac{11}{10}}\langle
2^k\delta\rangle^{\frac{4}{5}}\right).
\end{eqnarray*}
Using the Hardy-Littlewood-Sobolev, H\"{o}lder's inequality,
$(\ref{e75})$, we estimate ${\rm II_{12}}$ as the following :
\begin{eqnarray*}
 {\rm II_{12}} & \leq
& \|\chi_k u\|_{L_t^2L_x^5}\big\||\nabla|^{-\frac{3}{2}}(u u_{> \frac{M}{16}})\big\|_{L_t^\infty L_x^2} \\
& \lesssim &
\|\chi_k|y|^{-2}\|_{L_t^4L_x^5}\big\||y|^2u\big\|_{L_t^4L_x^\infty}
\|u u_{>\frac{M}{16}}\|_{L_t^\infty L_x^{5/4}}\\
& \lesssim_u & M^{-1}(2^k\delta)^{-\frac{3}{4}}\langle 2^k\delta
\rangle^{\frac{1}{2}} G\big(\frac{M}{16}\big).
\end{eqnarray*}
Thus, the left hand side of $(\ref{e714})$ can be bounded by:
\[\big(N^{-\frac{6}{5}}\delta^{-\frac{11}{10}} +
N^{-\frac{6}{5}}\delta^{-\frac{3}{10}} +
 N^{-2}\delta^{-\frac{3}{2}}+ N^{-1}\delta^{-\frac{3}{4}} +
 N^{-1}\delta^{-\frac{1}{2}}\big)G\big(\frac{N}{16}\big).\]
This is acceptable by choosing $N$ sufficiently large depending on
$\delta$ and $\eta$.

The last claim follows from the time reversal symmetry and the
$L_x^2$-boundedness of $P^{\pm}$.

We now turn to the region of $(t,y)$ integration where $|y| \ll
M|t|$. To begin with, we recall the bounds in \cite{c15} for the
kernels of the propagators in the region  $|x| \leq N^{-1}$,  $|y|
\ll M|t|$,  $|t| \geq \delta \gg N^{-2}$; and the region $|x| \geq
N^{-1}$, $y$ and $t$ as above:
\begin{eqnarray*}
 \big|P_Me^{-it\Delta}(x,y)\big| + \big|P_M^{\pm}e^{-it\Delta}(x,y)\big| \lesssim
 \frac{1}{(M^2|t|)^{50}}K_M(x,y),
\end{eqnarray*}
where\[K_M(x,y):= \dfrac{M^5}{\langle M(x-y)\rangle^{50}} +
\dfrac{M^5}{\langle Mx\rangle^2\langle My\rangle^2\langle
M|x|-M|y|\rangle^{50}}\] be bounded on $L_x^2$.

Let $\tilde \chi_k$ be  the characteristic function of the set
\[\{\,(t,y): 2^k\delta \leq |t| \leq 2^{k+1}\delta,\, |y| \ll M|t|\,\}.\]
\begin{lemma}[The tail]
Let $\eta > 0$ be a small number and $\delta$ be as in Lemma
$8.3$. Then
\[\sum_{M\geq N}\sum_{k=0}^\infty
\Big\|\int_{2^k\delta}^{2^{k+1}\delta}\int_{\mathbb R^5}
\frac{K_M(x,y)}{(M^2|t|)^{50}}\tilde\chi_k(t,y)[\tilde P_M
|\nabla|^{\frac{1}{2}}F(u(t))](y)\,\mathrm{d} y \,\mathrm{d} t
\Big\|_{L_x^2} \lesssim_u \eta G\big(\frac{N}{16}\big)\] for
sufficiently large $N$ depending on $u$ and $\eta$.
\end{lemma}
{\it Proof}. By Minkowski's inequality, the boundedness of $K_M$,
the support property of $\tilde P_M$, H\"{o}lder's and the
Hardy-Littlewood-Sobolev inequality
\begin{eqnarray*}
& & \Big\|\int_{2^k\delta}^{2^{k+1}\delta}\int_{\mathbb R^5}
\frac{K_M(x,y)}{(M^2|t|)^{50}}\tilde\chi_k(t,y)[\tilde P_M
|\nabla|^{\frac{1}{2}}F(u(t))](y)\,\mathrm{d} y \,\mathrm{d} t \Big\|_{L_x^2}\nonumber\\
& \lesssim & (M^22^k\delta)^{-50}\big\|\tilde \chi_k(t,y)[\tilde
P_M
|\nabla|^{\frac{1}{2}}F(u)]\big\|_{L_t^1L_y^2}\nonumber\\
& \lesssim & (M^22^k\delta)^{-50}2^k\delta
M^{\frac{1}{2}}\Big\|\tilde P_M \big(|\nabla|^{-2}\big( \,
\big|u_{\leq \frac{M}{16}} + u_{>\frac{M}{16}}\big|^2 \,
\big)(u_{\leq \frac{M}{16}}+ u_{>
\frac{M}{16}})\big)\Big\|_{L_t^\infty
L_y^2} \\
& \lesssim & (M^22^k\delta)^{-50}2^k\delta M^{\frac{1}{2}}\Big(
\big\| |\nabla|^{-2}(u^2)u_{> \frac{M}{16}}\big\|_{L_t^\infty
L_y^2} + \big\| |\nabla|^{-2}(u u_{>
\frac{M}{16}})u_{\leq \frac{M}{16}}\big\|_{L_t^\infty L_y^2}\Big)\\
& \lesssim & (M^22^k\delta)^{-50}2^k\delta M^{\frac{1}{2}}\Big(
\big\||\nabla|^{-2}(u^2)\big\|_{L_t^\infty
L_x^{5/2}}\|u_{>\frac{M}{16}}\|_{L_t^\infty L_x^{10}} \\
& & \hspace{4cm} + \big\||\nabla|^{-2}( u
u_{>\frac{M}{16}})\big\|_{L_t^\infty
L_x^{10}}\|u\|_{L_t^\infty L_x^{5/2}}\Big)\\
&\lesssim_u & (M^22^k\delta)^{-50}2^k\delta
M^{\frac{1}{2}}\Big(\|u\|^2_{L_t^\infty
L_x^{5/2}}M^{\frac{3}{2}}G\big(\frac{M}{16}\big) +
\|u_{>\frac{M}{16}}\|_{L_t^\infty
L_x^{10}}\|u\|_{L_t^\infty L_x^{5/2}}\Big)\\
& \lesssim_u & (M^22^k\delta)^{-50}2^k\delta M^2
G\big(\frac{M}{16}\big)
\end{eqnarray*}
Summing first over $k \geq 0$ and then $M \geq N$,
we obtain
\begin{eqnarray*}
\sum_{M\geq N}\sum_{k=0}^\infty
\Big\|\int_{2^k\delta}^{2^{k+1}\delta}\int_{\mathbb
R^5}\frac{K_M(x,y)}{(M^2|t|)^{50}} \tilde\chi_k(t,y)[\tilde P_M
|\nabla|^{\frac{1}{2}}F(u(t))](y)\,\mathrm{d} y \,\mathrm{d} t \Big\|_{L_x^2} \\
 \lesssim_u (N^2\delta)^{-49}G\big(\frac{N}{16}\big).\end{eqnarray*}
Choosing $N$ sufficiently large depending on $\delta, \eta$, we
get the desired result.

From $(\ref{e710})$, $(\ref{e712})$, Lemma 8.3, Lemma 8.4, Lemma
8.5, it concludes Proposition 8.1, which in turn proves Theorem
8.1.

\section{No soliton}
\setcounter{equation}{0}

In this section we prove
\begin{theorem}
There exists no non-zero soliton-like solution in the sense of
Theorem $1.3$.
\end{theorem}
{\it Proof}. We argue by contradiction. Assume that there exists
such a soliton solution, then by Theorem 6.1, Theorem 8.1, $u \in
L_t^\infty H^s_x(s \geq 1)$, and $u$ has the energy of the form
\[E(u(t)) = \frac{1}{2}\int_{\mathbb R^5} |\nabla u|^2 \, \mathrm{d} x -
 \frac{1}{4}\iint_{\mathbb R^5 \times \mathbb R^5}\frac{|u(x)|^2|u(y)|^2}{|x-y|^3} \, \mathrm{d}x \mathrm{d}y.\]

Now, define
\[M_a(t) :=  2 {\rm Im}\int_{\mathbb R^5} \bar u(t,x)\vec a(x) \cdot \nabla u(t,x) \,\mathrm{d} x,\]
where $a(x) = x\psi\big(\frac{|x|}{R}\big)$, $\psi$ is a smooth,
radial function such that
\[\psi(r)=
\begin{cases}
1, & r \leq 1\\
0, &  r \geq 2.
\end{cases} \]
 Then, by the Cauchy-Schwarz inequality, we have
\begin{equation}\label{e81}
|M_a(t)| \leq R\|u\|_2\|\nabla u\|_2 \lesssim_u R.
\end{equation}

We should prove by our assumption $\sup_{t \in \mathbb
R}\big\||\nabla|^{\frac{1}{2}}u\big\|_2 < \frac{\sqrt{6}}{3}\big\|
|\nabla|^{\frac{1}{2}}Q\big\|_2$ that $M_a(t)$ is an increasing
function of time, i.e., $\partial_t M_a(t) > 0$. Thus, a
contradiction with $(\ref{e81})$

A few computations with equation $(1.1)$ yields
\begin{eqnarray}
\partial_t M_a(t)& = & 12E(u(t)) - 2\int_{\mathbb R^5} |\nabla u|^2 \, \mathrm {d}x \label{e80}\\
& &  -\int_{\mathbb
R^5}\Big[\frac{24}{R|x|}\psi'\big(\frac{|x|}{R}\big) +
\frac{11}{R^2}\psi^{''}\big(\frac{|x|}{R}\big) +
\frac{|x|}{R^3}\psi^{'''}\big(\frac{|x|}{R}\big)\Big]|u(t,x)|^2 \,\mathrm{d} x \label{e82}\\
& & + 4\int_{\mathbb R^5}\Big[\psi\big(\frac{|x|}{R}\big) -1 +
\frac{|x|}{R}\psi'\big(\frac{|x|}{R}\big)\Big]|\nabla u(t,x)|^2 \,\mathrm{d} x \label{e83}\\
& & \hspace{-2.4cm} -3 \iint_{\mathbb R^5 \times \mathbb
R^5}\Big[x\psi\big(\frac{|x|}{R}\big) -
y\psi\big(\frac{|y|}{R}\big)-(x-y)\Big]\cdot
\frac{x-y}{|x-y|^5}|u(t,x)|^2|u(t,y)|^2 \,\mathrm{d} x
\,\mathrm{d} y \label{e84}.
\end{eqnarray}
We will prove that $(\ref{e82})$, $(\ref{e83})$, $(\ref{e84})$ are
sufficiently small compared to $(\ref{e80})$.

Note that $(\ref{e82})$ has a trivial bound $R^{-2}$.

Now, let $\eta > 0$ be a small number to be chosen later. From
Lemma 8.1, there exists $R = R(\eta)$ such that for all $t \in
\mathbb R$
\begin{equation}\label{e85}
\int_{|x| \geq \frac{R}{4}}|u(t,x)|^2 \,\mathrm{d} x \leq \eta.
\end{equation}
Define $\chi$ as a smooth cutoff to the region $|x| \geq
\frac{R}{2}$ with $\nabla \chi$ be bounded by $R^{-1}$ and
supported on $\{|x| \thicksim R\}$. Since $u \in C_t^0 H^s(s >
1)$, using the interpolation theorem and $(\ref{e85})$, we deduce
\begin{eqnarray*}|(\ref{e83})| \lesssim \| \chi\nabla u(t)\|^2_2
\lesssim  & \|\nabla(\chi u)\|_2^2 + \|u \nabla \chi\|_2^2
\lesssim
\|\chi u(t)\|_2^{\frac{2(s-1)}{s}}\|u(t)\|_{H^s}^{\frac{2}{s}} + \eta\\
 \lesssim_u \eta^{\frac{s-1}{s}} +\eta.
\end{eqnarray*}
It remains to estimate $(\ref{e84})$. We divide the integration
into three parts.
\[\aligned
(\ref{e84}) & =   \\
 & 2\mu \int\!\!\!\int_{\substack{ |x|\geq R \\  |y|\geq R
}}
\bigg(x\Big(\psi\big(\frac{|x|}{R}\big)-1\Big)-y\Big(\psi\big(\frac{|y|}{R}\big)-1\Big)\bigg)\cdot
\frac{x-y}{|x-y|^5}|u(t,x)|^2|u(t,y)|^2 \,\mathrm{d} x \,\mathrm{d} y\\
& + 2\mu\iint_{\substack{ |x|\geq R \\  |y|< R }}
x\Big(\psi\big(\frac{|x|}{R}\big)-1\Big)\cdot
\frac{x-y}{|x-y|^5}|u(t,x)|^2|u(t,y)|^2 \,\mathrm{d} x \,\mathrm{d} y\\
& - 2\mu \iint_{\substack{ |x|< R \\  |y|\geq R }}
y\Big(\psi\big(\frac{|y|}{R}\big)-1\Big)\cdot
\frac{x-y}{|x-y|^5}|u(t,x)|^2|u(t,y)|^2 \,\mathrm{d} x \,\mathrm{d} y\\
& := I_1 + I_2 + I_3.
\endaligned\]
We first estimate $I_1$. By the Gagliardo-Nirenberg inequality of
convolution type and $(\ref{e85})$
\[|I_1| \lesssim \iint_{\substack{ |x|\geq R \\  |y|\geq R
}}\frac{|u(x)|^2|u(y)|^2}{|x-y|^3}\,\mathrm{d} x \,\mathrm{d} y
\lesssim \|\chi u\|_2\|\nabla u\|_2^3 \lesssim_u \eta^{1/2}.\] To
estimate $I_2$, using the Hardy-Littlewood-Sobolev  inequality,
Lemma 3.1, Sobolev's embedding theorem,
\begin{eqnarray*} |I_2| & \lesssim &
\iint_{\substack{ |x|>2R \\ |y|< R }}
|x|\frac{|u(x)|^2|u(y)|^2}{|x-y|^4} \,\mathrm{d} x \,\mathrm{d} y \\
& &+ \iint_{\substack{ R< |x| \leq 2R \\ |y|< R }}
\bigg|x\Big(\psi\big(\frac{|x|}{R}\big)-1\Big)\bigg|
\frac{|u(x)|^2|u(y)|^2}{|x-y|^4} \,\mathrm{d} x \,\mathrm{d} y\\
& \lesssim & \iint_{\mathbb R^5 \times \mathbb R^5}\frac{|\chi
u(x)|^2|u(y)|^2}{|x-y|^3} \,\mathrm{d} x \,\mathrm{d} y +
R^{-\frac{3}{4}}\iint_{\mathbb R^5 \times
\mathbb R^5}\frac{|x|^{7/4}|u|\cdot|\chi u(x)||u(y)|^2}{|x-y|^4} \,\mathrm{d} x \,\mathrm{d} y\\
& \lesssim & \|\chi u\|_2\|\nabla(\chi u)\|_2\|\nabla u\|_2^2 +
R^{-\frac{3}{4}}\big\||x|^{7/4}u\big\|_{L^\infty_x}\|\chi u\|_2\| u\|_{H^1_x}^2 \\
& \hspace{6pt}\lesssim_u &  \eta^{\frac{2s-1}{2s}} +
R^{-\frac{3}{4}}\eta^{\frac{1}{2}}.
\end{eqnarray*}
Note that in the last inequality, we used the interpolation as
that  to estimate $(\ref{e83})$.

$I_3$ can be estimated in the same argument.

Thus, choosing $\eta$ sufficiently small depending on $u$, $R$
sufficiently large depending on $u$ and $\eta$, we have
$$|(\ref{e82})| + |(\ref{e83})| + |(\ref{e84})| \lesssim
\frac{1}{100} \times \bigg[12E(u(t)) - 2\int_{\mathbb R^5} |\nabla
u|^2 \, \mathrm {d}x\bigg].$$

On the other hand, as $\sup_{t \in \mathbb
R}\big\||\nabla|^{\frac{1}{2}}u\big\|_2 < \frac{\sqrt 6}{3}
\big\||\nabla|^{\frac{1}{2}}Q\big\|_2$, using the
Hardy-Littlewood-Sobolev type inequality $(\ref{a1})$, we see
$(\ref{e80}) > 0$. Hence $\partial_t M_a(t) > 0$.  This concludes
the proof of Theorem 9.1.

\vskip0.5cm

\textbf{Acknowledgements:} The authors would like to thank Prof. B.
Pausader  for his invaluable comments and suggestions.   The authors
are partly supported by the NSF of China (No. 10725102 and No.
10726053).

\end{document}